\documentclass[a4paper,12pt,oneside]{article}
\usepackage[english]{babel}
\usepackage[T1]{fontenc} 
\usepackage[utf8]{inputenc}
\usepackage{amsthm}
\usepackage{amsmath}
\usepackage{amssymb}  
\usepackage{indentfirst}
\usepackage{fancyhdr}
\usepackage{amsthm}
\usepackage{graphicx}
\usepackage{pdfpages}
\usepackage{esint}
\usepackage{url}

\definecolor{dg}{rgb}{0.01, 0.75, 0.24}

\numberwithin{equation}{section}

\theoremstyle{plain} 
\newtheorem{thm}{Theorem}[section]

\newtheorem{lem}[thm]{Lemma} 
\newtheorem{prop}[thm]{Proposition}

\def\io{\int_{\Omega}}
\def\ibtt{\int_{B_{t'}}}
\def\ibt {\int_{B_{t}}}

\def\ibr2 {\int_{B_{\frac{R}{2}}}}
\def\ibro2 {\int_{B_{\frac{\rho}{2}}}}

\def\fhi{\varphi}
\def\R{\mathbb{R}}
\def\N{\mathbb{N}}
\def\A{\mathcal{A}}
\def\dd{\textrm{d}}

\def\fibr2 {\displaystyle\fint_{B_{\frac{R}{2}}}}
\def\fibro2 {\fint_{B_{\frac{\rho}{2}}}}

\usepackage{geometry}
\allowdisplaybreaks[4]

\begin{document}

\title{\textbf{Higher differentiability for bounded solutions to a class of obstacle problems with $(p,q)$-growth}}

\author{
\sc{Antonio Giuseppe Grimaldi}\thanks{
 Dipartimento di Matematica e Applicazioni "R. Caccioppoli", Università degli Studi di Napoli "Federico II", Via Cintia, 80126 Napoli
 (Italy). E-mail: \textit{antoniogiuseppe.grimaldi@unina.it}}}

\maketitle
\maketitle

\begin{abstract}
We establish the higher fractional differentiability of bounded minimizers to a class of obstacle problems with non-standard growth conditions of the form
\begin{gather*}
\min \biggl\{ \displaystyle\int_{\Omega} F(x,Dw)dx \ : \ w \in \mathcal{K}_{\psi}(\Omega)  \biggr\},
\end{gather*}
where $\Omega$ is a bounded open set of $\mathbb{R}^n$, $n \geq 2$,
the function $\psi \in W^{1,p}(\Omega)$ is a fixed function called \textit{obstacle} and 
$
\mathcal{K}_{\psi}(\Omega) := \{ w \in W^{1,p}(\Omega) : w \geq \psi \ \text{a.e. in} \ \Omega  \}
$
is the class of admissible functions. If the obstacle $\psi$ is locally bounded, we prove that the gradient of solution inherits some fractional differentiability property, assuming that both the gradient of the obstacle and the mapping $x \mapsto D_\xi F(x,\xi)$ belong to some suitable Besov space.
The main novelty is that such assumptions are not related to the dimension $n$.
\end{abstract}


\section{Introduction}
In this paper we study the higher fractional differentiability properties of the gradient of bounded solutions $u \in W^{1,p}(\Omega)$ to obstacle problems of the form
\begin{gather}\label{!1}
\min \biggl\{ \displaystyle\int_{\Omega} F(x,Dw)dx \ : \ w \in \mathcal{K}_{\psi}(\Omega)  \biggr\}.
\end{gather}
Here $\Omega$ is a bounded open set of $\mathbb{R}^n$, $n \geq 2$,
the function $\psi: \Omega \rightarrow [-\infty, + \infty)$, called \textit{obstacle}, belongs to the Sobolev class $W^{1,p}(\Omega)$ and 
\begin{gather}\label{!2}
\mathcal{K}_{\psi}(\Omega) := \{ w \in W^{1,p}(\Omega) : w \geq \psi \ \text{a.e. in} \ \Omega  \}
\end{gather}
is the class of admissible functions.
Note that the set $\mathcal{K}_{\psi}(\Omega) $ is not empty since $\psi \in \mathcal{K}_{\psi}(\Omega) $.

In what follows, we assume that the energy density $F : \Omega \times \mathbb{R}^n \rightarrow [0, + \infty)$ is a Carath\'{e}odory function with Uhlenbeck structure, i.e. there exists a function $\tilde{F} : \Omega \times [0, + \infty) \rightarrow [0, + \infty)$ such that
$$ F(x, \xi)= \tilde{F}(x, |\xi|)\eqno{\rm{{ (F1)}}}$$
for a.e. $x \in \Omega$ and every $\xi \in \mathbb{R}^n$.
\\Moreover, we also assume that there exist positive constants $\tilde{\nu}$, $\tilde{L}$, $\tilde{l}$, exponents $2 \leq p < q < + \infty$ and a parameter $\mu \in [0,1]$, that will allow us to consider in our analysis both the degenerate and the non-degenerate situation, such that the following assumptions are satisfied:
$$ \dfrac{1}{\tilde{l}} (|\xi|^p-\mu^p ) \leq F(x, \xi) \leq \tilde{l}(\mu^2+|\xi|^2)^{\frac{q}{2}}\eqno{\rm{{ (F2)}}}$$
$$  \langle D_{\xi \xi}F(x, \xi) \lambda, \lambda \rangle \geq \tilde{\nu} (\mu^2+ |\xi|^2)^{\frac{p-2}{2}} |\lambda|^2  \eqno{\rm{{ (F3)}}}$$
$$ |D_{\xi \xi }F(x, \xi)| \leq \tilde{L}(\mu^2 + |\xi|^2)^{\frac{q-2}{2}}\eqno{\rm{{ (F4)}}}$$
for a.e. $x,y \in \Omega$ and every $\xi \in \mathbb{R}^n$.
\\Very recently, in \cite{cupini.marcellini.mascolo.passarelli} it has been proved that, under the structure assumption (F1), (F3) and (F4) imply (F2), i.e. if $p<q$, the functional $F$ has non-standard growth conditions of $(p,q)$-type, as initially defined and studied by Marcellini \cite{marcellini1986,marcellini1991}.

We say that function $F$ satisfies assumption (F5) if there exists a non-negative function $g \in L^{\frac{p+2\beta}{p+\beta-q}}_{\text{loc}}(\Omega)$, with $0 < \beta< \alpha <1$, such that
$$ |D_{\xi}F(x, \xi)- D_{\xi}F(y, \xi)| \leq |x-y|^{\alpha}(g(x)+g(y))(\mu^2+ |\xi|^2)^{\frac{q-1}{2}}\eqno{\rm{{ (F5)}}}$$
for a.e. $x,y \in \Omega$ and every $\xi \in \mathbb{R}^n$.\\
On the other hand, we say that assumption (F6) is satisfied if there exists a sequence of measurable non-negative functions $g_k \in L^{\frac{p+2\alpha}{p+\alpha-q}}_{\text{loc}}(\Omega)$, with $0 <  \alpha <1$, such that
$$\displaystyle\sum_{k=1}^{\infty} \Vert g_k \Vert^{\sigma}_{L^{\frac{p+2\alpha}{p+\alpha-q}}(\Omega)} < \infty,$$
for some $\sigma \ge 1$,
and at the same time
$$|D_{\xi}F(x,\xi)-D_{\xi}F(y, \xi)| \leq |x-y|^{\alpha} (g_k(x)+g_k(y)) (\mu^2 +|\xi|^2)^{\frac{q-1}{2}} \eqno{\rm{{ (F6)}}} $$
for a.e. $x,y \in \Omega$ such that $2^{-k} \text{diam}(\Omega) \leq |x-y| < 2^{-k+1}\text{diam}(\Omega)$ and for every $\xi \in \mathbb{R}^n$.\\
It is well known that the regularity of the minima often comes from the fact that they are also solutions to the corresponding Euler-Lagrange system, in the unconstrained setting, or, in the case of obstacle problems, to the corresponding variational inequality.
In the case of standard growth conditions, that is $p=q$, $u \in W^{1,p}(\Omega)$ is a solution to the obstacle problem in $\mathcal{K}_{\psi}(\Omega) $ if, and only if, $u \in \mathcal{K}_{\psi}(\Omega) $ solves the variational inequality 
\begin{align}\label{1}
\displaystyle\int_{\Omega} \langle \mathcal{A}(x,Du)  ,D(\varphi-u) \rangle dx \geq 0 
\end{align}
for all $\varphi \in \mathcal{K}_{\psi}(\Omega)$, where we set
\begin{gather}
\mathcal{A}(x,\xi)= D_{\xi}F(x,\xi). \notag
\end{gather}
On the other hand, in the case of non-standard growth conditions, even in the unconstrained case, the relation between extremals and minima is an issue that requires a careful investigation
(see for example \cite{carozza1,carozza2,eleuteri.passarelli2021}).

From conditions (F2)--(F4), we deduce the existence of positive constants $\nu,L,l$ such that the following $p$-ellipticity and $q$-growth conditions are satisfied by the map $\mathcal{A}$:
$$|\mathcal{A}(x,\xi)| \leq l(\mu^2 +|\xi|^2)^{\frac{q-1}{2}} \eqno{\rm{{ (A1)}}}$$
$$\langle \mathcal{A}(x,\xi)-\mathcal{A}(x,\eta), \xi-\eta \rangle \geq \nu |\xi-\eta|^{2} (\mu^2+|\xi|^{2}+|\eta|^{2})^{\frac{p-2}{2}} \eqno{\rm{{ (A2)}}}$$
$$|\mathcal{A}(x,\xi)-\mathcal{A}(x,\eta)|\leq L|\xi-\eta| (\mu^2+|\xi|^{2}+|\eta|^{2})^{\frac{q-2}{2}} \eqno{\rm{{ (A3)}}}$$
for a.e. $x,y \in \Omega$, for every $\xi, \eta \in \mathbb{R}^n$. 
\\Furthermore, if condition (F5) 
or (F6)
holds, then $\A$ satisfies assumption (A4),
or (A5), respectively
i.e.
$$|\mathcal{A}(x,\xi)-\mathcal{A}(y, \xi)| \leq |x-y|^{\alpha} (g(x)+g(y)) (\mu^2 +|\xi|^2)^{\frac{q-1}{2}} \eqno{\rm{{ (A4)}}}$$
for a.e. $x,y \in \Omega$ and every $\xi \in \mathbb{R}^n$.
or
$$|\A(x,\xi)-\A(y, \xi)| \leq |x-y|^{\alpha} (g_k(x)+g_k(y)) (\mu^2 +|\xi|^2)^{\frac{q-1}{2}} \eqno{\rm{{ (A5)}}}$$
for a.e. $x,y \in \Omega$ such that $2^{-k} \text{diam}(\Omega) \leq |x-y| < 2^{-k+1}\text{diam}(\Omega)$ and for every $\xi \in \mathbb{R}^n$.

The study of obstacle problems started with the works by Stampacchia \cite{stampacchia} and Fichera \cite{fichera} and has since then attracted much attention.
It is usually observed that the regularity of the solutions to the obstacle problems is influenced by the one of the obstacle; for example, for linear obstacle problems, obstacle and solutions have the same regularity \cite{brezis.kinderlehrer,caffarelli.kinderlehrer,kinderlehrer.stampacchia}. This does not apply in the nonlinear setting, hence there have been intense research activities in this direction (see \cite{choe,choe.l,eleuteri.h,michael}, just to mention a few).

In recent years, there has been a considerable interest in analyzing how an extra differentiability of integer or fractional order of the obstacle transfers to the gradient of solutions: for instance we quote \cite{caselli.gentile.giova,eleuteri.passarelli,eleuteri.passarelli1,gentile.0,gentile0,grimaldi0} in the setting of standard growth conditions, \cite{defilippis1,defilippis,gavioli1,gavioli2,gentile1,grimaldi.ipocoana,grimaldi.ipocoana1,zhang.zheng} in the setting of non-standard growth conditions.

The analysis comes from the fact that the regularity of the solutions to obstacle problem \eqref{1} is strictly connected to the analysis of the regularity of the solutions to partial differential equations of the form
$$\text{div} \mathcal{A}(x,Du)= \text{div} \mathcal{A}(x,D\psi),$$
whose higher differentiability properties have been widely investigated (see for instance \cite{pambrosio,baison.clop2017,clop,giova01,giova,giova1,passarelli2014.potanal,passarelli2014}.\\
It is well known that no extra differentiability properties for the solutions can be expected even if the obstacle $\psi$ is smooth, unless some assumption is given on the coefficients of the operator $\mathcal{A}$.
Therefore, recent results concerning the higher differentiability of solutions to obstacle problems show that a $W^{1,r}$ Sobolev regularity, with $r \ge n$, or a $B^s_{r,\sigma}$ Besov regularity, with $r \ge \frac{n}{s}$, on the partial map $x \mapsto \mathcal{A}(x,\xi)$ is a sufficient condition (see \cite{eleuteri.passarelli,gavioli1,gavioli2,gentile.0} for the case of Sobolev class of integer order and \cite{eleuteri.passarelli,grimaldi.ipocoana} for the fractional one).

When referring to functionals with non-standard growth conditions, by looking at the counterexamples in \cite{giaquinta,marcellini1991}, if the ratio $q/p \nrightarrow 1$ when $n \rightarrow \infty$, then minimizers become unbounded.
On the other hand, many regularity results require that $q/p \rightarrow 1$ when $n \rightarrow \infty$. Now, it is well known that, both for unconstrained and constrained problems with $(p,q)$-growth, when dealing with bounded minimizers, regularity results for the gradient can be proved under dimension-free conditions on the gap $q/p$ and weaker assumptions on the data of the problem (see \cite{carozza0,colombo1,gentile0,gentile1,giova1}). Moreover, in
\cite{gentile0,gentile1,giova1}, the higher differentiability of integer order of bounded solutions to \eqref{1} is obtained assuming that the coefficients of $\mathcal{A}$ and the gradient of the obstacle belong to a Sobolev class that is not related to the dimension $n$ but to the ellipticity and the growth exponents of the functional.

Recently, in \cite{derosa.grimaldi}, it has been proved that, assuming the local boundedness of the obstacle $\psi $, the solution to obstacle problem \eqref{!1} is locally bounded under a sharp relation between $p$ and $q$.
Here, we study higher fractional differentiability properties of bounded solutions to obstacle problems satisfying $(p,q)$-growth conditions. The novelty of this work consists in showing that, even in the fractional setting, the higher differentiability properties of bounded solutions to \eqref{!1} hold true assuming that the Besov type regularity on the partial map $x \mapsto \mathcal{A}(x,\xi)$  and on the gradient of the obstacle are not related to the dimension $n$. We observe that the bound \eqref{gap2} is only needed to get the local boundedness of the solution (see Theorem \ref{thmbound}). Therefore, if we deal with a priori bounded minimizers, then the result holds without the hypotesis \eqref{gap2}.

More precisely, we shall prove the following theorems.

\begin{thm}\label{mainthm1}
Let $F(x,\xi)$ satisfy (F1)--(F5) for exponents $2\leq p  <q$ such that
\begin{equation}\label{gap}
q <  p+\beta
\end{equation}
and
\begin{equation}\label{gap2}
  	\frac{1}{q}\ge \frac{1}{p}-\frac{1}{n-1}.
  \end{equation}
Let $u\in \mathcal{K}_\psi(\Omega)$ be the solution to the obstacle problem \eqref{!1}. If $\psi \in L^\infty_{loc}(\Omega)$, then we have
\begin{equation}
D\psi \in B^\alpha_{\frac{p+2\beta}{p+1+\beta-q},\infty,\textrm{loc}}(\Omega)\Rightarrow (\mu^2 +|Du|^2)^\frac{p-2}{4} Du \in B^\alpha_{2,\infty,\textrm{loc}}(\Omega),
\end{equation}
provided $0< \beta< \alpha < 1$.
\end{thm}

On the other hand, a Besov regularity of the type $B^\alpha_{r,\sigma}$, with $\sigma $ finite,
is stronger than the one of the type $B^\alpha_{p,\infty}$. In this case, we prove higher fractional differentiability properties for bounded minimizers under weaker assumptions both on the coefficients of $\mathcal{A}$ and on the gradient of the obstacle and on the bound for the gap $q/p$. The main difference is that a stronger embedding theorem between Sobolev and Besov spaces holds (see Lemma \ref{lemmasobbes}).

\begin{thm}\label{mainthm2}
Let $F(x,\xi)$ satisfy (F1)--(F4) and (F6) for exponents $2\leq p  <q$ verifying \eqref{gap2} and
\begin{equation}\label{gap3}
q <  p+ \alpha .
\end{equation}
Let $u\in \mathcal{K}_\psi(\Omega)$ be the solution to the obstacle problem \eqref{!1}. If $\psi \in L^\infty_{loc}(\Omega)$, then we have
\begin{equation}
D\psi \in B^\alpha_{\frac{p+2 \alpha }{p+1+ \alpha-q},\sigma,\textrm{loc}}(\Omega)\Rightarrow (\mu^2 +|Du|^2)^\frac{p-2}{4} Du \in B^{ \alpha}_{2,\sigma,\textrm{loc}}(\Omega),
\end{equation}
provided $\sigma(1+\alpha) \le 2$.
\end{thm}

The structure of this paper is the following. After recalling some notation and preliminary results in Section \ref{secnot}, 
a Gagliardo-Niremberg type inequality in Besov spaces is established in Section \ref{iiB} for a priori bounded minimizers. Then,  
we concentrate on proving our main result. The strategy is to prove uniform a priori estimates for solutions to a family of approximating problems. 
Therefore, in Section \ref{appl}, we present the approximation lemma that allows to construct a sequence of functions satisfying $p$-growth conditions that converges to the initial problem.
In Section \ref{pthm1}, we prove Theorem \ref{mainthm1}. In particular, we derive the a priori estimates in Section \ref{apriori}, and, in Section \ref{limit}, we pass to the limit in the approximating problems. Eventually, in Section \ref{teorema1.2}, we give the proof of Theorem \ref{mainthm2}, focusing only on the a priori estimate since the approximation procedure works exactly in same way as in the proof of Theorem \ref{mainthm1}.

The local boundedness allows us to use an interpolation inequality (see Lemma \ref{interineq1}) that gives the higher local integrability $L^{p+2\alpha}$ of the gradient of the solutions. Such higher integrability is the key tool in order to weaken the assumptions on the function $g $ and on $D \psi$ with respect to the higher differentiability result established in \cite{grimaldi.ipocoana}. Indeed, for $p< n-2\alpha$ and $q< p+\alpha-\frac{\alpha(p+2\alpha)}{n}$, we have $  L^\frac{n}{\alpha} \subset L^\frac{p+2\alpha}{p+\alpha-q}  $, and, moreover, under our assumption on the gap, i.e. $q< p +\alpha$, $  B^\alpha_{2q-p,\sigma} \subset B^\alpha_{\frac{p+2\alpha}{p+1+\alpha-q},\sigma} $.

\section{Notations and preliminary results}
\label{secnot}
For the rest of the paper, we denote by $C$, $c$, $\pi$ general positive constants. Different occurrences from line to line will be still denoted using the same letters. Relevant dependencies on parameters will be emphasized using parentheses or subscripts. 
We denote by $B(x,r)=B_{r}(x)= \{ y \in \mathbb{R}^{n} : |y-x | < r  \}$ the ball centered at $x$ of radius $r$. We shall omit the dependence on the center and on the radius when no confusion arises. For a function $u \in L^{1}(B)$, the symbol
\begin{center}
$\displaystyle\fint_{B} u(x) dx = \dfrac{1}{|B|} \displaystyle\int_{B} u(x) dx$.
\end{center}
will denote the integral mean of the function $u$ over the set $B$.

It is convenient to introduce an auxiliary function
\begin{center}
$V_{p}(\xi)=(\mu^{2}+|\xi|^{2})^\frac{p-2}{4} \xi$
\end{center}
defined for all $\xi\in \mathbb{R}^{n}$. One can easily check that, for $p \geq 2$, it holds
\begin{gather}
|\xi|^p \leq  |V_p(\xi)|^2. \label{Vp}
\end{gather}
For the auxiliary function $V_{p}$, we recall the following estimate (see the proof of \cite[Lemma 8.3]{giusti}): 
\begin{lem}\label{D1}
Let $1<p<+\infty$. There exists a constant $c=c(n,p)>0$ such that
\begin{center}
$c^{-1}(\mu^{2}+|\xi|^{2}+|\eta|^{2})^{\frac{p-2}{2}} \leq \dfrac{|V_{p}(\xi)-V_{p}(\eta)|^{2}}{|\xi-\eta|^{2}} \leq c(\mu^{2}+|\xi|^{2}+|\eta|^{2})^{\frac{p-2}{2}} $
\end{center}
for any $\xi, \eta \in \mathbb{R}^{n}$, $\xi \neq \eta$.
\end{lem}

Now we state a well-known iteration lemma (see \cite[Lemma 6.1]{giusti} for the proof).
\begin{lem}\label{lm2}
Let $\Phi  :  [\frac{R}{2},R] \rightarrow \mathbb{R}$ be a bounded nonnegative function, where $R>0$. Assume that for all $\frac{R}{2} \leq r < s \leq R$ it holds
$$\Phi (r) \leq \theta \Phi(s) +A + \dfrac{B}{(s-r)^2}+ \dfrac{C}{(s-r)^{\gamma}}$$
where $\theta \in (0,1)$, $A$, $B$, $C \geq 0$ and $\gamma >0$ are constants. Then there exists a constant $c=c(\theta, \gamma)$ such that
$$\Phi \biggl(\dfrac{R}{2} \biggr) \leq c \biggl( A+ \dfrac{B}{R^2}+ \dfrac{C}{R^{\gamma}}  \biggr).$$
\end{lem}

The following regularity result, whose proof can be found in \cite{derosa.grimaldi} in a more general setting, allows us to obtain the local boundedness of solutions to obstacle problem \eqref{!1}.

\begin{thm}\label{thmbound}
Let $u \in W^{1,p}(\Omega)$ be a solution to \eqref{!1} under assumptions (F2) and (F3), for exponents $2 \leq p < q $ verifying \eqref{gap2}.
If $\psi \in L^\infty_{\text{loc}}(\Omega)$, then $u \in L^\infty_{\text{loc}}(\Omega)$ and the following estimate
\begin{equation}\label{limiteq}
    \sup_{B_{R/2}}|u| \leq C (\sup_{B_{R}}|\psi|+\Vert u \Vert_{W^{1,p}(B_{R})})^\pi,
\end{equation}
holds for every ball $B_{R} \Subset \Omega$, for $\pi:= \pi(n,p,q)$ and with $C:=C(n,p,q,R)$.
\end{thm}

\subsection{Difference quotient}
\label{secquo}
We recall some properties of the finite difference quotient operator that will be needed in the sequel. Let us recall that, for every function $F:\mathbb{R}^{n}\rightarrow \mathbb{R}$ the $r$-th finite difference operator is defined by
\begin{align*}
\tau^1_{s,h}F(x):= & \tau_{s,h}F(x) =F(x+he_{s})-F(x),\\
\tau^r_{s,h}F(x):= & \tau_{s,h} (\tau^{r-1}_{s,h}F(x)), \quad r \in \N, r \geq 1,
\end{align*}
where $h \in \mathbb{R}^{n}$, $e_{s}$ is the unit vector in the $x_{s}$ direction and $s \in \{1,...,n\}$.
\\We start with the description of some elementary properties that can be found, for example, in \cite{giusti}.
\begin{prop}\label{rapportoincrementale}
Let $F \in W^{1,p}(\Omega)$, with $p \geq1$, and let us consider the set
\begin{center}
$\Omega_{|h|} = \{ x \in \Omega : \text{dist}(x,\partial \Omega)> |h|  \}$.
\end{center}
Then
\\(i) $\tau_{h}F \in W^{1,p}(\Omega_{|h|})$ and 
\begin{center}
$D_{i}(\tau_{h}F)=\tau_{h}(D_{i}F)$.
\end{center}
(ii) If at least one of the functions $F$ or $G$ has support contained in $\Omega_{|h|}$, then
\begin{center}
$$\displaystyle\int_{\Omega}F \tau_h G dx = \displaystyle\int_{\Omega} G \tau_{-h}F dx.$$
\end{center}
(iii) We have $$\tau_h (FG)(x)= F(x+h)\tau_h G(x)+G(x) \tau_h F(x).$$
\end{prop}
The next result about finite difference operator is a kind of integral version of Lagrange Theorem.
\begin{lem}\label{ldiff}
If $0<\rho<R,$ $|h|<\frac{R-\rho}{2},$ $1<p<+\infty$ and $F,\ DF \in L^{p}(B_{R})$, then
\begin{center}
$\displaystyle\int_{B_{\rho}} |\tau_{h}F(x)|^{p} dx \leq c(n,p)|h|^{p} \displaystyle\int_{B_{R}} |DF(x)|^{p}dx$.
\end{center}
Moreover,
\begin{center}
$\displaystyle\int_{B_{\rho}} |F(x+h)|^{p} dx \leq  \displaystyle\int_{B_{R}} |F(x)|^{p}dx$.
\end{center}
\end{lem}

\subsection{Besov-Lipschitz spaces}
\label{secbesov}
Let $v:\mathbb{R}^{n} \rightarrow \mathbb{R}$ be a function. As in \cite[Section 2.5.12]{haroske}, given $\alpha >0 $ and $1 \leq p,s< \infty$, we say that $v$ belongs to the Besov space $B^{\alpha}_{p,s}(\mathbb{R}^{n})$ if $v \in L^{p}(\mathbb{R}^{n})$ and
\begin{center}
$\Vert v \Vert_{B^{\alpha}_{p,s}(\mathbb{R}^{n})} = \Vert v \Vert_{L^{p}(\mathbb{R}^{n})} + [v]_{B^{\alpha}_{p,s}(\mathbb{R}^{n})} < \infty$,
\end{center}
where
\begin{center}
$[v]_{B^{\alpha}_{p,s}(\mathbb{R}^{n})} =  \biggl( \displaystyle\int_{\mathbb{R}^{n}} \biggl( \displaystyle\int_{\mathbb{R}^{n}} \dfrac{|\tau^r_hv(x)|^{p}}{|h|^{\alpha p}} dx \biggr)^{\frac{s}{p}}  \dfrac{dh}{|h|^{n}} \biggr)^{\frac{1}{s}}  < \infty$.
\end{center}
Here and in what follows, $r$ is the smallest integer larger than $\alpha$.
Equivalently, we could simply say that $v \in L^{p}(\mathbb{R}^{n})$ and $\frac{\tau^r_{h}{v}}{|h|^{\alpha}} \in L^{s}\bigl( \frac{dh}{|h|^{n}}; L^{p}(\mathbb{R}^{n}) \bigr)$. As usual, if one integrates for $h \in B(0, \delta)$ for a fixed $\delta >0$ then an equivalent norm is obtained, because
\begin{center}
$\biggl( \displaystyle\int_{\{|h| \geq \delta\}} \biggl( \displaystyle\int_{\mathbb{R}^{n}} \dfrac{|\tau^r_h v(x)|^{p}}{|h|^{\alpha p}} dx \biggr)^{\frac{s}{p}}  \dfrac{dh}{|h|^{n}} \biggr)^{\frac{1}{s}} \leq c(n, \alpha,p,s, \delta) \Vert v \Vert_{L^{p}(\mathbb{R}^{n})} $.
\end{center}
Similarly, we say that $v \in B^{\alpha}_{p,\infty}(\mathbb{R}^{n})$ if $v \in L^{p}(\mathbb{R}^{n})$ and
\begin{center}
$[v]_{B^{\alpha}_{p, \infty}(\mathbb{R}^{n})} =  \displaystyle\sup_{h \in \mathbb{R}^{n}} \biggl( \displaystyle\int_{\mathbb{R}^{n}} \dfrac{|\tau^r_hv(x)|^{p}}{|h|^{\alpha p}} dx \biggr)^{\frac{1}{p}} < \infty $.
\end{center}
Again, one can simply take supremum over $|h| \leq \delta$ and obtain an equivalent norm. By construction, one has $B^{\alpha}_{p, s}(\mathbb{R}^{n}) \subset L^{p}(\mathbb{R}^{n})$. One also has the following version of Sobolev embeddings (a proof can be found at \cite[Proposition 7.12]{haroske}).
\begin{lem}\label{3.1}
Suppose that $0 < \alpha <1$.
\\ (a) If $1 < p < \frac{n}{\alpha}$ and $1 \leq s \leq p^{*}_{\alpha} = \frac{np}{n- \alpha p}$, then there is a continuous embedding $B^{\alpha}_{p, s}(\mathbb{R}^{n}) \subset L^{p^{*}_{\alpha}}(\mathbb{R}^{n})$.
\\ (b) If $p = \frac{n}{\alpha}$ and $1 \leq s \leq \infty$, then there is a continuous embedding $B^{\alpha}_{p, s}(\mathbb{R}^{n}) \subset BMO(\mathbb{R}^{n})$,
\\ where $BMO$ denotes the space of functions with bounded mean oscillations \emph{\cite[Chapter 2]{giusti}}.
\end{lem}
For further needs, we recall the following inclusions (\cite[Proposition 7.10 and Formula (7.35)]{haroske}).
\begin{lem}\label{3.2}
Suppose that $0 < \beta < \alpha $.
\\ (a) If $1 < p < \infty$ and $1 \leq s \leq t \leq \infty$, then $B^{\alpha}_{p, s}(\mathbb{R}^{n}) \subset B^{\alpha}_{p, t}(\mathbb{R}^{n})$.
\\ (b) If $1 < p < \infty$ and $1 \leq s , t \leq \infty$, then $B^{\alpha}_{p, s}(\mathbb{R}^{n}) \subset B^{\beta}_{p, t}(\mathbb{R}^{n})$.
\end{lem}

Given a domain $\Omega \subset \mathbb{R}^{n}$, we say that $v$ belongs to the local Besov space $ B^{\alpha}_{p, s,loc}$ if $\varphi \ v \in B^{\alpha}_{p, s}(\mathbb{R}^{n})$ whenever $\varphi \in \mathcal{C}^{\infty}_{c}(\Omega)$. It is worth noticing that one can prove suitable version of Lemma \ref{3.1} and Lemma \ref{3.2}, by using local Besov spaces.

We have the following lemma, which can be found in \cite{baison.clop2017} in the case $\alpha \in (0,1)$.
\begin{lem}
A function $v \in L^{p}_{loc}(\Omega)$ belongs to the local Besov space $B^{\alpha}_{p,s,loc}$ , with $0< \alpha < 2$, if, and only if,
\begin{center}
$\biggl\Vert \dfrac{\tau^r_{h}v}{|h|^{\alpha}} \biggr\Vert_{L^{s}\bigl(\frac{dh}{|h|^{n}};L^{p}(B)\bigr)}<  \infty,$
\end{center}
for any ball $B\subset2B\subset\Omega$ with radius $r_{B}$. Here the measure $\frac{dh}{|h|^n}$ is restricted to the ball $B(0,r_B)$ on the h-space.
\end{lem}
\proof 
Let $\alpha \in [1, 2)$. Using twice Proposition \ref{rapportoincrementale} (iii), for any smooth and compactly supported test function $\varphi$, the following pointwise identity
\begin{align}
    \frac{\tau^2_h(\varphi v)(x)}{|h|^\alpha} & =
    \frac{\tau_h(\varphi(x+h) \tau_hv(x)+v(x) \tau_h\varphi(x))}{|h|^\alpha} \notag\\
    &=
\varphi(x+2h)\frac{\tau^2_hv(x)}{|h|^\alpha}+
2 \frac{\tau_hv(x)\tau_h\varphi(x+h)}{|h|^\alpha}+
v(x) \frac{\tau^2_{h}\varphi(x)}{|h|^\alpha} \label{pointwiseid}
\end{align} 
holds.
\\It is clear that
\begin{align}
    \left\vert\frac{\tau_hv(x)\tau_h\varphi(x+h)}{|h|^\alpha} \right\vert
    \leq &  \frac{|\tau_hv(x)|}{|h|^{\alpha-1} } \Vert  D \varphi \Vert_{\infty},
\end{align}
and so, since $\alpha-1 < 1$ and $v \in B^{\alpha-1}_{p,s,\text{loc}}(\Omega)$ from Lemma \ref{3.2} (b), one has 
$$\frac{\tau_hv(x)\tau_h\varphi(x+h)}{|h|^\alpha} \in L^{s}\biggl(\frac{dh}{|h|^{n}};L^{p}(\R^n)\biggr).$$
Moreover, it holds
\begin{align}
   \left\vert \frac{v(x)\tau^2_h\varphi(x)}{|h|^\alpha} \right\vert \leq & |v(x)| \Vert \tau_h D \varphi \Vert_{\infty}|h|^{1-\alpha} \notag\\
    \leq & |v(x)| \Vert  D^2 \varphi \Vert_{\infty}|h|^{2-\alpha}
\end{align}
and therefore we have
$$\frac{v(x)\tau^2_h\varphi(x)}{|h|^\alpha} \in L^{s}\biggl(\frac{dh}{|h|^{n}};L^{p}(\R^n)\biggr).$$
As a consequence, we have the equivalence
$$ \varphi v \in B^{\alpha}_{p,s}(\R^n) \Longleftrightarrow \varphi(x+2h)\frac{\tau^2_hv(x)}{|h|^\alpha} \in L^{s}\biggl(\frac{dh}{|h|^{n}};L^{p}(\R^n)\biggr). $$
However, it is clear that $\varphi(x+2h)\frac{\tau^2_hv(x)}{|h|^\alpha} \in L^{s}\biggl(\frac{dh}{|h|^{n}};L^{p}(\R^n)\biggr)$ for every $\varphi \in \mathcal{C}^\infty_c(\Omega)$ if, and only if, the same happens for every $\varphi = \chi_B$ and all ball $B \subset 2B \subset \Omega$. This concludes the proof.
\endproof

It is known that Besov-Lipschitz spaces of fractional order $\alpha \in (0,1)$ can be characterized in pointwise terms. Given a measurable function $v:\mathbb{R}^{n} \rightarrow \mathbb{R}$, a \textit{fractional $\alpha$-Hajlasz gradient for $v$} is a sequence $\{g_{k}\}_{k}$ of measurable, non-negative functions $g_{k}:\mathbb{R}^{n} \rightarrow \mathbb{R}$, together with a null set $N\subset\mathbb{R}^{n}$, such that the inequality 
\begin{center}
$|v(x)-v(y)|\leq (g_{k}(x)+g_{k}(y))|x-y|^{\alpha}$
\end{center} 
holds whenever $k \in \mathbb{Z}$ and $x,y \in \mathbb{R}^{n}\setminus N$ are such that $2^{-k} \leq|x-y|<2^{-k+1}$. We say that $\{g_{k}\}_{k} \in l^{s}(\mathbb{Z};L^{p}(\mathbb{R}^{n}))$ if
\begin{center}
$\Vert \{g_{k}\}_{k} \Vert_{l^{s}(L^{p})}=\biggl(\displaystyle\sum_{k \in \mathbb{Z}}\Vert g_{k} \Vert^{s}_{L^{p}(\mathbb{R}^{n})} \biggr)^{\frac{1}{s}}<  \infty.$
\end{center} 

The following result was proved in \cite{koskela}.
\begin{thm}
Let $0< \alpha <1,$ $1 \leq p < \infty$ and $1\leq s \leq \infty $. Let $v \in L^{p}(\mathbb{R}^{n})$. One has $v \in B^{\alpha}_{p,s}(\mathbb{R}^{n})$ if, and only if, there exists a fractional $\alpha$-Hajlasz gradient $\{g_{k}\}_{k} \in l^{s}(\mathbb{Z};L^{p}(\mathbb{R}^{n}))$ for $v$. Moreover,
\begin{center}
$\Vert v \Vert_{B^{\alpha}_{p,s}(\mathbb{R}^{n})}\simeq \inf \Vert \{g_{k}\}_{k} \Vert_{l^{s}(L^{p})},$
\end{center}
where the infimum runs over all possible fractional $\alpha$-Hajlasz gradients for $v$.
\end{thm}

\section{Gagliardo-Niremberg inequality}\label{iiB}
In this section, we collect some results in Besov spaces that will be useful later. 
\begin{lem}\label{lemmaint1}
Let $v \in W^{1,p}_{loc}(\R^n)$. If $Dv \in B^\gamma_{p,s,loc}(\R^n)$, for some $1 \le s \le \infty$ and $0< \gamma <1$, then $v \in B^{1+\gamma}_{p,s,loc}(\R^n)$. Moreover, the following estimate 
$$[v]_{B^{1+\gamma}_{p,s}(B_\rho)} \leq c [Dv]_{B^{\gamma}_{p,s}(B_R)}$$
holds for every ball $B_\rho \subset B_R $, with $c:=c(n,p)$.
\end{lem}
\proof We give the proof of Lemma \ref{lemmaint1} only for $s=\infty$, since the case $s $ finite can be obtained in a similar way.\\
Fix $0<\rho<R$, $|h|<\frac{R-\rho}{2}$ and consider balls $B_\rho \subset B_R$.
Since $1 < 1+\gamma < 2$, we have that
\begin{equation}
[v]_{B^{1+\gamma}_{p,\infty}(B_\rho)} =  \displaystyle\sup_{h \in \mathbb{R}^{n}} \biggl( \displaystyle\int_{B_\rho} \dfrac{|\tau^2_hv(x)|^{p}}{|h|^{(1+\gamma) p}} dx \biggr)^{\frac{1}{p}}. \notag
\end{equation}
Now, using the fact that $v \in W^{1,p}_{loc}(\R^n)$ and Lemma \ref{ldiff}, we obtain
\begin{align*}
\displaystyle\int_{B_\rho} \dfrac{|\tau_h(\tau_hv(x))|^{p}}{|h|^{(1+\gamma) p}} dx \leq & c |h|^p \displaystyle\int_{B_R} \dfrac{|\tau_hDv(x)|^{p}}{|h|^{(1+\gamma) p}} dx \\
= & c \displaystyle\int_{B_R} \dfrac{|\tau_hDv(x)|^{p}}{|h|^{\gamma p}} dx \le c [Dv]_{B^{\gamma}_{p,\infty}(B_R)}^p,
\end{align*}
which is finite by the assumption on $Dv$. This completes the proof.
\endproof

The following interpolation inequality can be found in \cite{runst}.

\begin{lem}\label{interineq1}
Let $1 \leq p < \infty$, $1 \leq s \leq \infty$, $\gamma >0$ and $0 < \theta < 1$. Then, the following interpolation inequality
\begin{equation}\label{eqinte1}
    \Vert  v \Vert_{B^{\theta \gamma}_{p/\theta,s/\theta}(\R^n)} \leq c  \Vert  v \Vert^\theta_{B^{ \gamma}_{p,s}(\R^n)} \Vert v  \Vert^{1-\theta}_{L^\infty(\R^n)} 
\end{equation}
holds for every $v \in B^{ \gamma}_{p,s}(\R^n) \cap L^\infty(\R^n)$.
\end{lem}

For the next result, see e.g. \cite{miro}.

\begin{lem}\label{besintosob}
Let $\gamma >0$, $1 \leq p < \infty$ and $1 \leq s \leq \infty$. If $m \in \N$, $\gamma > m$, then there is a continuous embedding $B^\gamma_{p,s}(\R^n) \subset W^{m,p}(\R^n)$.
\end{lem}

Now, we are able to prove the following higher integrability result.
\begin{prop}\label{lemma2.14}
Let $v \in W^{1,p}_{loc}(\R^n) \cap L^\infty_{loc}(\R^n)$ and let $Dv \in B^\gamma_{p,\infty,loc}(\R^n)$, for some 
$1 \leq p < \infty$ and $0< \gamma <1$. Then $Dv \in L^{p(1+\beta)}_{loc}(\R^n)$, for every $0 < \beta < \gamma$. Moreover, the following estimate
\begin{align*}
    \int_{B_\rho} |Dv|^{p(1+\beta)}dx 
    \le &  C \Vert v \Vert^{p\beta}_{L^\infty(B_R)}
    \biggl(
    [ Dv ]^{p}_{B^{\gamma}_{p,\infty}(B_R)}
     + 
    \frac{1}{(R-\rho)^{2p}}
    \Vert v \Vert^{p}_{W^{1,p}(B_R)}
    \biggr)
\end{align*}
holds for every ball $B_\rho \subset B_{R} $, with $C:=C(n,p,\gamma, \beta)$.
\end{prop}
\proof
Thanks to Lemma \ref{lemmaint1}, we obtain
$$ v \in B^{1+ \gamma}_{p,\infty}(\R^n) \quad \text{locally}.$$
Then, Lemma \ref{interineq1} yields
\begin{equation}
 v \in B^{\theta(1+ \gamma)}_{p/\theta,\infty}(\R^n) \quad \text{locally}, 
\end{equation}
for every $\theta \in (0,1)$.\\
Choosing $\theta = \frac{1}{1+ \beta}$, for $0 <\beta < \gamma$, we have
$$\theta (1+\gamma)= \frac{1+\gamma}{1+ \beta}>1.$$
Let us consider $0 < \rho < R \le 1$ and fix balls $B_\rho \subset B_{R}$ and a cut-off function $\eta \in \mathcal{C}_c^\infty(B_\frac{R+\rho}{2})$, $\eta =1$ on $B_\rho$ such that $|D\eta| \le \frac{C}{R-\rho}$ and $|D^2\eta| \le \frac{C}{(R-\rho)^2}$. By virtue of Lemma \ref{besintosob}, we have
\begin{equation}
    \int_{B_\rho} |Dv|^{p(1+\beta)} dx \le \Vert \eta v \Vert^{p(1+\beta)}_{W^{1,p(1+\beta)}(\R^n)}
    \le c  \Vert \eta v \Vert^{p(1+\beta)}_{B^{\frac{1+\gamma}{1+\beta}}_{p(1+\beta),\infty}(\R^n)}
    .\label{embedding}
\end{equation}
From Lemma \ref{interineq1}, we get
\begin{align}
     \Vert \eta v \Vert^{p(1+\beta)}_{B^{\frac{1+\gamma}{1+\beta}}_{p(1+\beta),\infty}(\R^n)}
     \le c \Vert v \Vert^{p\beta}_{L^\infty(B_R)}
      \Vert \eta v \Vert^p_{B^{1+\gamma}_{p,\infty}(\R^n)}. \label{embeddings1}
\end{align}
Using identity \eqref{pointwiseid} and properties of $\eta$, we infer
\begin{align}
   \Vert \eta v \Vert^{p}_{B^{1+\gamma}_{p,\infty}(\R^n)} 
   \le & C \Vert v \Vert^{p}_{L^{p}(B_\frac{R+\rho}{2})}
   +C \sup_{ |h| \le \frac{R-\rho}{4}} \int_{\R^n} |\eta(x+2h)|^p \frac{|\tau^2_hv|^p}{|h|^{p(1+\gamma)}} dx \notag\\
   &+ C \sup_{|h| \le \frac{R-\rho}{4}} \int_{\R^n}
   |\tau_h \eta(x+h)|^{p}
    \frac{|\tau_h v |^{p}}{|h|^{p(1+\gamma)}} dx
    + C \sup_{|h| \le \frac{R-\rho}{4}} 
    \int_{B_\frac{3R+\rho}{4}} |v|^p \frac{|\tau^2_h\eta(x)|^p}{|h|^{p(1+\gamma)}} dx \notag\\
    \le & C \Vert v \Vert^{p}_{L^{p}(B_\frac{R+\rho}{2})}
   +C \sup_{ |h| \le \frac{R-\rho}{4}} \int_{B_\frac{3R+\rho}{4}}  \frac{|\tau^2_hv|^p}{|h|^{p(1+\gamma)}} dx \notag\\
   &+ C \sup_{|h| \le \frac{R-\rho}{4}} \int_{\R^n}
   \Vert D \eta \Vert^p_{L^\infty(B_\frac{3R+\rho}{4})}
    \frac{|\tau_h v |^{p}}{|h|^{p\gamma}} dx \notag\\
    &
    + C \sup_{|h| \le \frac{R-\rho}{4}} |h|^{p(1-\gamma)}
    \int_{B_\frac{3R+\rho}{4}} |v|^p \Vert D^2 \eta\Vert^p_{L^\infty(B_\frac{3R+\rho}{4})}  dx \notag\\
   \le &
   C \Vert v \Vert^{p}_{L^{p}(B_\frac{R+\rho}{2})}
   + C [  v ]^{p}_{B^{1+\gamma}_{p,\infty}(B_\frac{3R+\rho}{4})} \notag\\
   &
    +\frac{C}{(R-\rho)^{p}} \sup_{|h| \le \frac{R-\rho}{4}} \int_{B_\frac{3R+\rho}{4}}
    \frac{|\tau_h v |^{p}}{|h|^{p\gamma}} dx \notag\\
    &
    +\frac{C}{(R-\rho)^{2p}} \sup_{|h| \le \frac{R-\rho}{4}}
    |h|^{p(1 -\gamma)}\int_{B_\frac{3R+\rho}{4}} |v |^{p} dx. \notag
    \end{align}
Now, exploiting Lemma \ref{ldiff} and using the fact that $R-   \rho <1$, we obtain    
    \begin{align}
     \Vert \eta v \Vert^{p(1+\beta)}_{B^{\frac{1+\gamma}{1+\beta}}_{p(1+\beta),\infty}(\R^n)}
    \le &
      C \Vert v \Vert^{p}_{L^{p}(B_\frac{R+\rho}{2})}
   + C [  v ]^{p}_{B^{1+\gamma}_{p,\infty}(B_\frac{3R+\rho}{4})} \notag\\
   &
    +\frac{C}{(R-\rho)^{p}} \sup_{|h| \le \frac{R-\rho}{4}} |h|^{p(1-\gamma)} \int_{B_R}
    |D v |^{p} dx \notag\\
    &
    +\frac{C}{(R-\rho)^{2p}}\int_{B_\frac{3R+\rho}{4}} |v |^{p} dx \notag\\
    \le & 
    C [ v ]^{p}_{B^{1+\gamma}_{p,\infty}(B_{\frac{3R+\rho}{4}})} +
    \frac{C}{(R-\rho)^{2p}} \Vert v \Vert^{p}_{W^{1,p}(B_R)}.
    \label{lastinequality}
\end{align}
Combining inequalities \eqref{embedding}, \eqref{embeddings1} and \eqref{lastinequality} and Lemma \ref{lemmaint1}, we derive

\begin{align*}
    \int_{B_\rho} |Dv|^{p(1+\beta)}dx \leq &  
     C \Vert v \Vert^{p\beta}_{L^\infty(B_R)}
     [ v ]^{p}_{B^{1+\gamma}_{p,\infty}(B_{\frac{3R+\rho}{4}})}
     + 
    \frac{C}{(R-\rho)^{2p}}
    \Vert v \Vert^{p\beta}_{L^\infty(B_R)}
    \Vert v \Vert^{p}_{W^{1,p}(B_R)} \notag\\
    \le &  C \Vert v \Vert^{p\beta}_{L^\infty(B_R)}
    [ Dv ]^{p}_{B^{\gamma}_{p,\infty}(B_R)}
     + 
    \frac{C}{(R-\rho)^{2p}}
    \Vert v \Vert^{p\beta}_{L^\infty(B_R)}
    \Vert v \Vert^{p}_{W^{1,p}(B_R)}
\end{align*}
i.e. the desired estimate.
\endproof

The following proposition is an immediate consequence of the previous result and will be a fundamental tool for the proof of Theorem \ref{mainthm1}.
\begin{prop}\label{mainint}
Let $v \in W^{1,p}_{loc}(\R^n) \cap L^\infty_{loc}(\R^n)$, for some $p \geq 2$, and assume that $V_p(Dv) \in B^\gamma_{2,\infty,loc}(\R^n)$, for some $0< \gamma <1$. Then $Dv \in L^{p+2 \beta}_{loc}(\R^n)$, for every $0 < \beta < \gamma$. Moreover, the following inequality
\begin{align*}
    \int_{B_\rho} |Dv|^{p+2\beta}dx 
    \le &  C \Vert v \Vert^{2\beta}_{L^\infty(B_R)}
    \biggl(
    [ V_p(Dv) ]^{2}_{B^{\gamma}_{2,\infty}(B_R)}
     + 
    \frac{1}{(R-\rho)^{2p}}
    \Vert v \Vert^{p}_{W^{1,p}(B_R)}
    \biggr)
\end{align*}
holds for every ball $B_\rho \subset B_{R} $, with $C:=C(n,p,\gamma,\beta)$.
\end{prop}
\proof
By Lemma \ref{D1}, we get
\begin{equation}\label{bes1}
    |\tau_hV_p(Dv)|^2 \geq c |\tau_hDv|^2 (\mu^2+|Dv(x+h)|^2+|Dv(x)|^2)^{\frac{p-2}{2}} \geq c |\tau_hDv|^p,
\end{equation}
where in the last inequality we used the fact that $p \geq 2$.\\
Estimate \eqref{bes1} implies 
\begin{equation}\label{bes2}
    \int_{B_R} \dfrac{|\tau_h Dv|^p}{|h|^{p \frac{2 \gamma}{p}}} = \int_{B_R} \dfrac{|\tau_h Dv|^p}{|h|^{2 \gamma}} dx \le c \int_{B_R} \dfrac{|\tau_h V_p(Dv)|^2}{|h|^{2 \gamma}} dx \quad \text{for every } h.
\end{equation}
Now, taking the supremum over $h$ in \eqref{bes2} and by virtue of the assumption on $V_p(Dv)$, we derive that
\begin{equation}\label{bes3}
    Dv \in B^{\frac{2 \gamma}{p}}_{p,\infty}(\R^n) \quad \text{locally}.
\end{equation}
Thanks to Lemma \ref{lemmaint1} and \eqref{bes3}, we obtain
$$ v \in B^{1+\frac{2 \gamma}{p}}_{p,\infty}(\R^n) \quad \text{locally}.$$
By virtue of Lemma \ref{lemma2.14} and \eqref{bes2}, it follows
\begin{align*}
    \int_{B_\rho} |Dv|^{p+2\beta}dx 
    \le &  C \Vert v \Vert^{2\beta}_{L^\infty(B_R)}
    [ V_p(Dv) ]^{2}_{B^{\gamma}_{2,\infty}(B_R)}
     + 
    \frac{C}{(R-\rho)^{2p}}
    \Vert v \Vert^{2\beta}_{L^\infty(B_R)}
    \Vert v \Vert^{p}_{W^{1,p}(B_R)}
\end{align*}
for every $0 < \beta < \gamma$.
\endproof

Moreover, we have the following embeddings  between Sobolev and Besov spaces (see \cite{trebel}).
\begin{lem}\label{lemmasobbes}
Let $\gamma >0$, $1 < p < \infty$ and $1 \leq s \leq \min\{p,2\}$. Then, there is a continuous embedding $B^\gamma_{p,s}(\R^n) \subset W^{\gamma,p}(\R^n)$.
\end{lem}

Arguing similarly, but assuming that $V_p(Dv) \in B^\gamma_{2,s,loc}(\R^n)$, for some $1 \le s < \infty$, we are able to prove the following result.

\begin{prop}\label{higherintfin}
Let $v \in W^{1,p}_{loc}(\R^n) \cap L^\infty_{loc}(\R^n)$, for some $p \geq 2$, and assume that $V_p(Dv) \in B^\gamma_{2,s,loc}(\R^n)$, for some $1 \le s(1+\gamma) \le 2$ and $0< \gamma <1$. Then $Dv \in L^{p+2 \gamma}_{loc}(\R^n)$. Moreover, the following inequality
\begin{align*}
    \int_{B_\rho} |Dv|^{p+2\gamma}dx 
    \le &  C \Vert v \Vert^{2\gamma}_{L^\infty(B_R)}
    \biggl(
    [ V_p(Dv) ]^{2}_{B^{\gamma}_{2,s}(B_R)}
     + 
    \frac{1}{(R-\rho)^{2p}}
    \Vert v \Vert^{p}_{W^{1,p}(B_R)}
    \biggr)
\end{align*}
holds for every ball $B_\rho \subset B_{R} $, with $C:=C(n,p,\gamma)$.
\end{prop}

\section{Approximation lemma}\label{appl}
The main tool to prove Theorem \ref{mainthm1} 
is the following approximation lemma (see \cite{grimaldi.ipocoana} for the proof).

\begin{lem}\label{apprlem1}
Let $F : \Omega \times \R^n \rightarrow [0,+\infty), F = F(x, \xi)$, be a Carath\'{e}odory function satisfying assumptions (F1), (F2), (F3) and (F5). Then there exists
a sequence $(F_j)$ of Carath\'{e}odory functions $F_j: \Omega \times \R^n \rightarrow [0,+\infty)$ monotonically convergent to $F$, such that
\begin{itemize}
\item[(i)]for a.e. $x \in \Omega$ and every $\xi \in \R^n$, 
$F_j(x,\xi) = \tilde{F}_j (x,|\xi|)$,
\item[(ii)]for a.e. $x \in \Omega$, for every 	$\xi \in \R^n$ and for every $j$, 
$F_j(x,\xi)\leq F_{j+1}(x, \xi) \leq F(x,\xi)$,
\item[(iii)]for a.e. $x \in \Omega$ and every $\xi \in \R^n$, we have
$\langle D_{\xi\xi} F_j(x,\xi)\lambda, \lambda\rangle \geq \bar{\nu}(\mu^2+|\xi|^2)^\frac{p-2}{2}|\lambda|^2$,
with $\bar{\nu}$ depending only on $p$ and $\nu$,
\item[(iv)] for a.e. $x\in \Omega$ and for every $\xi \in \R^n$, there exist $L_1$, independent of $j$, and $\bar{L}_1$, depending on $j$, such that
\begin{align*}
& 1/L_1(|\xi|^p-\mu^p) \leq F_j(x,\xi)\leq L_1(\mu + |\xi|)^q,\\
&F_j (x,\xi)\leq \bar{L}_1(j)(\mu + |\xi|)^p,
\end{align*} 
\item[(v)]there exists a constant $C(j) > 0$ such that
\begin{align*}
&|D_\xi F_j(x,\xi)-D_\xi F_j(y,\xi)| \leq |x-y|^{\alpha} (k(x)+k(y))(\mu^2 +|\xi|^2)^\frac{q-1}{2},\\
&|D_\xi F_j(x,\xi)-D_\xi F_j(y,\xi)| \leq  C(j)|x-y|^{\alpha} (k(x)+k(y))(\mu^2 +|\xi|^2)^\frac{p-1}{2}
\end{align*}
for a.e. $x,y \in \Omega$ and for every $\xi \in \R^n$.
\end{itemize}
Assuming (F6) instead of (F5), statement $(v)$ would change as follows.
\begin{itemize}
\item[(v)]There exists a constant $C(j) > 0$ such that
\begin{align*}
&|D_\xi F_j(x,\xi)-D_\xi F_j(y,\xi)|\leq |x-y|^{\alpha} (g_k(x)+g_k(y))(\mu^2 +|\xi|^2)^\frac{q-1}{2},\\
&|D_\xi F_j(x,\xi)-D_\xi F_j(y,\xi)| \leq  C(j)|x-y|^{\alpha} (g_k(x)+g_k(y))(\mu^2 +|\xi|^2)^\frac{p-1}{2}
\end{align*}
for a.e. $x,y \in \Omega$ such that $2^{-k} \text{diam}(\Omega) \leq |x-y| < 2^{-k+1}\text{diam}(\Omega)$ and for every $\xi \in \R^n$.
\end{itemize}
\end{lem}

\section{Proof of Theorem \ref{mainthm1}}\label{pthm1}
This section is devoted to the proof of Theorem \ref{mainthm1}. In particular, in Section \ref{apriori}, we derive the a priori estimates for regular minimizers of obstacle problems \eqref{!1}, while in Section \ref{limit}, we conclude through an approximation argument.

\subsection{A priori estimate}\label{apriori}

We have the following theorem.

\begin{thm}\label{approximation}
Let $F(x,\xi)$ satisfy (F1)--(F5) for exponents $2\leq p < q$ such that \eqref{gap} and \eqref{gap2} hold.
Let $u\in \mathcal{K}_{\psi}(\Omega)$ be the solution to the obstacle problem \eqref{!1}.
Suppose that $$g \in L^{\frac{p+2\beta}{p+\beta-q}}_{loc}(\Omega)  , \ \psi \in L^\infty_{loc}(\Omega) \quad \text{and} \quad D \psi \in B^{\alpha}_{\frac{p+2\beta}{p+1+\beta-q}, \infty,loc}(\Omega),$$ for $0< \beta < \alpha  <1$. If we a priori assume that $$V_p(Du) \in B^\alpha_{2,\infty,loc}(\Omega),$$ then the following estimates
\begin{align}
    \int_{B_{R/4}} |Du|^{p+2\beta} dx 
     \leq &
     C (\Vert\psi\Vert_{L^\infty(B_R)}+\Vert u \Vert_{W^{1,p}(B_{R})})^\pi \notag\\
     & \cdot
     \biggl(
     \int_{B_{R}} (g^{\frac{p+2\beta}{p+\beta-q}}+1) dx +\Vert D\psi \Vert_{B^\alpha_{\frac{p+2\beta}{p+1+\beta-q},\infty}(B_R)}
     \biggr)^\pi  \label{stimapriori..}
\end{align}
and
\begin{align}
\displaystyle\int_{B_{R/4}} |\tau_h V_p(Du)|^2 dx 
     \le &  C  |h|^{2 \alpha} (\Vert\psi\Vert_{L^\infty(B_R)}+\Vert u \Vert_{W^{1,p}(B_{R})})^\pi \notag\\
     & \cdot
     \biggl(
     \int_{B_{R}} (g^{\frac{p+2\beta}{p+\beta-q}}+1) dx +\Vert D\psi \Vert_{B^\alpha_{\frac{p+2\beta}{p+1+\beta-q},\infty}(B_R)}
     \biggr)^\pi \label{stimapriori.}
\end{align}
hold for all balls $B_{R/4} \subset B_R \Subset \Omega$, for positive constants $C := C(n,p,q, \nu, L,R)$ and $\pi:= \pi (n,p,q,\beta)$.
\end{thm}\label{ape}
\proof By virtue of assumption \eqref{gap2} and Theorem \ref{thmbound}, $u \in L^\infty_{loc}(\Omega)$.
Hence, using Proposition \ref{mainint}, we deduce that
\begin{equation}\label{maxint}
    Du \in L^{p+2 \beta}_{loc}(\Omega).
\end{equation}
Notice that $Du \in L^{p+2 \beta}_{loc}(\Omega)$ implies that the $u$ satisfies the variational inequality \eqref{1} for every $\varphi \in W^{1,q}(\Omega)$ such that $\fhi \ge \psi$. Indeed, let $\varphi \in W^{1,q}(\Omega)$, $\fhi \ge \psi$, then the function
$ u+\varepsilon(\fhi -u)$ belongs to the admissible class, for every $\varepsilon \in (0,1)$, since 
\begin{align*}
u+\varepsilon(\fhi -v) = \varepsilon \fhi +(1-\varepsilon) u \geq \psi.
\end{align*}
Hence, by minimality of $u$, we get
\begin{align*}
\int_\Omega F(x,Du) d x \leq \int_\Omega F(x,Du+\varepsilon D(\fhi-u)) d x,
\end{align*}
which leads to
\begin{align*}
\int_\Omega [ F(x,Du+\varepsilon D(\fhi-u)) -F(x,Du) ] d x \geq 0.
\end{align*}
From Lagrange's theorem, for $\theta \in (0,1)$ it holds
\begin{align*}
\int_\Omega \langle \mathcal{A}(x,Du+\varepsilon \theta D(\fhi-u)), \varepsilon D(\fhi-u)\rangle d x\geq 0.
\end{align*}
Since $\varepsilon > 0$, we get
\begin{align}\label{varineqfr1}
\int_\Omega \langle \mathcal{A}(x,Du+\varepsilon \theta D(\fhi-u)), D(\fhi-u)\rangle \dd x\geq 0.
\end{align}
Now, from assumption (A1), we obtain
\begin{align}\nonumber
|\langle & \mathcal{A}(x,Du+\varepsilon \theta D(\fhi-u)),D(\fhi-u)\rangle|\\ \nonumber
& \le |\mathcal{A}(x,Du+\varepsilon \theta D(\fhi-u))||D(\fhi-u)|\\ \nonumber
& \le C(1+|Du+\varepsilon \theta D(\fhi-u)|^{q-1})
|D(\fhi-u)| \\ \nonumber
& \le C(1+|Du|^q+|D\fhi|^q), \nonumber
\end{align}
where we also used that $\varepsilon, \theta \in (0,1)$. 
\\On the other hand, by virtue of assumption \eqref{gap} and \eqref{maxint}, we have
$$1+|Du|^q+|D\fhi|^q \in L^1_{loc}(\Omega).$$
Therefore, by applying the Dominated convergence theorem, we can pass to the limit for $\varepsilon \to 0^+$ in \eqref{varineqfr1}, getting the inequality \eqref{1}, for every $\varphi \in W^{1,q}(\Omega)$ such that $\fhi \ge \psi$.

Fix $0 < \frac{R}{4} < \rho  < s < t < t'  <  \frac{R}{2} $ such that $B_{R} \Subset \Omega$ and a cut-off function $\eta \in \mathcal{C}_0^1(B_t)$ such that $0 \leq \eta \leq 1$, $\eta =1$ on $B_s$, $|D \eta | \leq \frac{C}{t-s}$.
\\Now, for $|h| \leq t'-t$, we consider functions
\begin{gather}
v_{1}(x)= \eta^{2}(x) [(u-\psi)(x+h)-(u-\psi)(x)] \notag 
\end{gather}
and
\begin{gather}
v_{2}(x)= \eta^{2}(x-h) [(u-\psi)(x-h)-(u-\psi)(x)]. \notag 
\end{gather}
Then
\begin{equation}
\varphi_1(x)=u(x)+tv_1 (x), \label{2:2} 
\end{equation}
\begin{equation}
\varphi_2(x)=u(x)+tv_2(x) \label{2:3}
\end{equation}
are admissible test functions for all $t \in [0,1)$.
\\Inserting \eqref{2:2} and \eqref{2:3} in \eqref{1}, we obtain
\begin{align}
\displaystyle\int_{\Omega}  \langle \mathcal{A}(x,Du), & D(\eta^2 \tau_h(u- \psi)) \rangle dx + \displaystyle\int_{\Omega} \langle \mathcal{A}(x,Du), D(\eta^2(x-h) \tau_{-h}(u- \psi)) \rangle dx \geq 0 .\label{2:4}
\end{align}
By means of a simple change of variable, we can write the second integral on the left hand side of the previous inequality as follows
\begin{align}
\displaystyle\int_{\Omega} \langle \mathcal{A}(x+h,Du(x+h)), D(-\eta^2 \tau_h(u- \psi)) \rangle dx \label{2:5}
\end{align}
and so inequality \eqref{2:4} becomes
\begin{align}
\displaystyle\int_{\Omega} \langle \mathcal{A}(x+h,Du(x+h))-\mathcal{A}(x,Du(x)), D(\eta^2 \tau_h(u- \psi)) \rangle dx 
\leq 0 \label{2:6}
\end{align}
We can write previous inequality as follows
\begin{align}
0 \geq & \displaystyle\int_{\Omega} \langle \mathcal{A}(x+h,Du(x+h))-\mathcal{A}(x+h,Du(x)),\eta^{2}D\tau_{h}u \rangle dx\notag\\
  &-\displaystyle\int_{\Omega} \langle \mathcal{A}(x+h,Du(x+h))-\mathcal{A}(x+h,Du(x)),\eta^{2}D\tau_{h}\psi \rangle dx\notag\\
  &+\displaystyle\int_{\Omega} \langle \mathcal{A}(x+h,Du(x+h))-\mathcal{A}(x+h,Du(x)),2\eta D \eta\tau_{h}(u-\psi) \rangle dx\notag\\
  &+\displaystyle\int_{\Omega} \langle \mathcal{A}(x+h,Du(x))-\mathcal{A}(x,Du(x)),\eta^{2}D\tau_{h}u \rangle dx\notag\\
  &-\displaystyle\int_{\Omega} \langle \mathcal{A}(x+h,Du(x))-\mathcal{A}(x,Du(x)),\eta^{2}D\tau_{h}\psi \rangle dx\notag\\
  &+\displaystyle\int_{\Omega} \langle \mathcal{A}(x+h,Du(x))-\mathcal{A}(x,Du(x)),2\eta D \eta\tau_{h}(u-\psi) \rangle dx \notag\\
 =:& I_{1}+I_{2}+I_{3}+I_{4}+I_{5}+I_{6}, \label{2:7}
\end{align}
that yields
\begin{align}
I_1 \leq & |I_2| + |I_3| + |I_4| + |I_5| + |I_6| . \label{2:8}
\end{align}
The ellipticity assumption (A2) and Lemma \ref{D1} imply
\begin{align}
I_{1} \geq &\nu \displaystyle\int_{\Omega}  \eta^{2} |\tau_{h}Du|^{2}(\mu^2 + |Du(x+h)|^{2}+|Du(x)|^{2})^{\frac{p-2}{2}} dx \notag\\
\geq & C(\nu) \displaystyle\int_{\Omega}  \eta^{2} |\tau_{h}V_p(Du)|^{2} dx . \label{I1}
\end{align}
From the growth condition (A3), Young's and H\"{o}lder's inequalities, Lemma \ref{D1} and assumption on $D \psi$, we get
\begin{align}\label{I2}
|I_{2}|\leq & L \displaystyle\int_{\Omega} \eta^{2} |\tau_{h}Du|(\mu^2 + |Du(x+h)|^{2}+|Du(x)|^{2})^{\frac{q-2}{2}}|\tau_{h}D \psi| dx \notag \\
\leq & \varepsilon\displaystyle\int_{\Omega} \eta^2 |\tau_hDu|^2 (\mu^2+|Du(x+h)|^2+|Du(x)|^2)^{\frac{p-2}{2}}dx \notag\\
&+ C_{\varepsilon}(L) \displaystyle\int_{\Omega} \eta^2 |\tau_h D \psi|^2 (\mu^2+|Du(x+h)|^2+|Du(x)|^2)^{\frac{2q-p-2}{2}} dx\notag\\
\leq & \varepsilon \displaystyle\int_{\Omega} \eta^2 |\tau_hV_p(Du)|^2 dx \notag\\
&+ C_{\varepsilon}(L) \biggl( \displaystyle\int_{B_t}|\tau_h D \psi|^{\frac{p+2\beta}{p+1+\beta-q}}dx \biggr)^{\frac{2(p+1+\beta-q)}{p+2\beta}}  \biggl( \displaystyle\int_{B_{t'}}(1+|D u|)^{p+2\beta}dx \biggr)^{\frac{2q-p-2}{p+2\beta}}\notag\\
\leq & \varepsilon \displaystyle\int_{\Omega} \eta^2 |\tau_hV_p(Du)|^2 dx \notag\\
&+ C_{\varepsilon}(L)|h|^{2\alpha}[D \psi]^2_{B^{\alpha}_{\frac{p+2\beta}{p+1+\beta-q}, \infty}(B_R)} \biggl( \displaystyle\int_{B_{t'}}(1+|D u|)^{p+2\beta}dx \biggr)^{\frac{2q-p-2}{p+2\beta}}. 
\end{align}
Arguing analogously, we get
\begin{align}
|I_{3}|  \leq & 2L \displaystyle\int_{\Omega} |D \eta| \eta |\tau_{h} Du| (\mu^2 + |Du(x+h)|^{2}+|Du(x)|^{2})^{\frac{q-2}{2}}|\tau_{h}(u-\psi)| dx \notag\\
\leq & \varepsilon \displaystyle\int_{\Omega} \eta^2 |\tau_hDu|^2 (\mu^2+|Du(x+h)|^2+|Du(x)|^2)^{\frac{p-2}{2}}dx \notag\\
&+ \dfrac{C_{\varepsilon}(L)}{(t-s)^2} \displaystyle\int_{B_t} |\tau_h(u- \psi)|^2(\mu^2+ |Du(x+h)|^2+|Du(x)|^2)^{\frac{2q-p-2}{2}}dx \notag\\
\leq & \varepsilon \displaystyle\int_{\Omega} \eta^2 |\tau_hV_p(Du)|^2 dx \notag\\
&+\dfrac{C_{\varepsilon}(L)}{(t-s)^2} \biggl( \displaystyle\int_{B_{t'}}|\tau_h (u-\psi)|^{\frac{p+2\beta}{p+1+\beta-q}}dx \biggr)^{\frac{2(p+1+\beta-q)}{p+2\beta}}  \biggl( \displaystyle\int_{B_{t'}}(1+|D u|)^{p+2\beta}dx \biggr)^{\frac{2q-p-2}{p+2\beta}} . \notag
\end{align}
Using Lemma \ref{ldiff}, we obtain
\begin{align}\label{I3}
|I_{3}| \leq & \varepsilon \displaystyle\int_{\Omega} \eta^2 |\tau_hV_p(Du)|^2 dx \notag\\
&+\dfrac{C_{\varepsilon}(L)}{(t-s)^2} |h|^2 \biggl( \displaystyle\int_{B_t}|D (u-\psi)|^{\frac{p+2\beta}{p+1+\beta-q}}dx \biggr)^{\frac{2(p+1+\beta-q)}{p+2\beta}}  \biggl( \displaystyle\int_{B_{t'}}(1+|D u|)^{p+2\beta}dx \biggr)^{\frac{2q-p-2}{p+2\beta}}.
\end{align}
In order to estimate the integral $I_4$, we use assumption (A4), Young's and H\"{o}lder's inequalities and Lemma \ref{D1} as follows
\begin{align}
|I_4| \leq & \displaystyle\int_{\Omega} \eta^2 |\tau_h Du| |h|^{\alpha} (g(x+h)+g(x))(1+|Du(x)|)^{\frac{q-1}{2}}dx \notag\\
\leq & \varepsilon \displaystyle\int_{\Omega} \eta^2 |\tau_hDu|^2 (\mu^2+|Du(x+h)|^2+|Du(x)|^2)^{\frac{p-2}{2}}dx \notag\\
&+ C_{\varepsilon}|h|^{2 \alpha} \displaystyle\int_{B_t} (g(x+h)+g(x))^2(1+|Du|)^{2q-p}dx \notag\\
\leq & \varepsilon \displaystyle\int_{\Omega} \eta^2 |\tau_hV_p(Du)|^2 dx \notag\\
&+ C_{\varepsilon}|h|^{2 \alpha} \biggl( \displaystyle\int_{B_{R}}g^\frac{p+2\beta}{p+\beta-q} dx \biggr)^{\frac{2(p+\beta-q)}{p2\beta}} \biggl(  \displaystyle\int_{B_t}(1+|Du|)^{p+2\beta}dx\biggr)^{\frac{2q-p}{p+2\beta}}. \label{I4}
\end{align}
We now take care of $I_5$. Similarly as above, exploiting assumption (A4) and H\"{o}lder's inequality, we infer
\begin{align*}
|I_5| 
\leq & \io \eta^2 |\tau_h D\psi| |h|^\alpha \left( g(x+h)+g(x)\right) \left( 1+|Du|^2\right)^{\frac{q-1}{2}} dx\\
\leq & |h|^\alpha \left( \ibtt g^\frac{p+2\beta}{p+\beta-q} dx \right)^\frac{p+\beta-q}{p+2\beta} \left( \ibt |\tau_h D\psi|^\frac{p+2\beta}{q+\beta}(1+|Du|)^\frac{(q-1)(p+2\beta)}{q+\beta} dx \right)^\frac{q+\beta}{p+2\beta}\\
\leq & |h|^\alpha \left( \ibtt g^\frac{p+2\beta}{p+\beta-q} dx \right)^\frac{p+\beta-q}{p+2\beta}
\left( \ibt |\tau_h D \psi|^\frac{p+2\beta}{p+1+\beta-q} dx \right)^\frac{p+1+\beta-q}{p+2\beta} \notag\\
& \cdot
\left( \ibt(1+|Du|)^\frac{(q-1)(p+2\beta)}{2q-p-1} dx \right)^\frac{2q-p-1}{p+2\beta}.
\end{align*}
Now, we observe
\begin{align}\label{disI5}
\frac{(q-1)(p+2\beta)}{2q-p-1} < p+2\beta \Longleftrightarrow p < q .
\end{align}
Hence
\begin{align}\label{I5}
|I_5| \leq |h|^{2\alpha}  \left( \int_{B_R} g^\frac{p+2\beta}{p+\beta-q} dx \right)^\frac{p+\beta-q}{p+2\beta}
[D \psi]_{B^\alpha_{\frac{p+2\beta}{p+1+\beta-q}, \infty}(B_R)}
\left( \ibt(1+|Du|)^{p+2\beta} dx \right)^\frac{q-1}{p+2\beta}.
\end{align}
From assumption (A4), hypothesis $|D\eta| < \frac{C}{t-s}$ and H\"{o}lder's inequality, we infer the following estimate for $I_6$.

\begin{align*}
|I_6| \leq &\frac{C}{t-s} |h|^\alpha \ibt |\tau_h (u-\psi)| (g(x+h)+g(x))(1+|Du|^2)^\frac{q-1}{2} dx\\
\leq &\frac{C}{t-s} |h|^\alpha  \left( \ibtt g^\frac{p+2\beta}{p+\beta-q} dx \right)^\frac{p+\beta-q}{p+2\beta}
\left( \ibt |\tau_h (u- \psi)|^\frac{p+2\beta}{p+1+\beta-q} dx \right)^\frac{p+1+\beta-q}{p+2\beta} \\
\cdot &
\left( \ibt(1+|Du|)^\frac{(q-1)(p+2\beta)}{2q-p-1} dx \right)^\frac{2q-p-1}{p+2\beta}.
\end{align*}
Using once again H\"{o}lder's inequality, inequality \eqref{disI5} and Lemma \ref{ldiff}, we have
\begin{align}
|I_6| \leq &\frac{C}{t-s} |h|^{\alpha+1}  \left( \int_{B_R} g^\frac{p+2\beta}{p+\beta-q} dx \right)^\frac{p+\beta-q}{p+2\beta}
\left( \ibtt |D (u- \psi)|^\frac{p+2\beta}{p+1+\beta-q} dx \right)^\frac{p+1+\beta-q}{p+2\beta} \notag\\
\cdot &
\left( \ibt(1+|Du|)^{p+2\beta} dx \right)^\frac{q-1}{p+2\beta} . \label{I6}
\end{align}
Inserting estimates \eqref{I1}, \eqref{I2}, \eqref{I3}, \eqref{I4}, \eqref{I5} and \eqref{I6} in \eqref{2:8}, we infer

\begin{align}
     C & (\nu) \displaystyle\int_{\Omega}  \eta^{2} |\tau_{h}V_p(Du)|^{2} dx \notag\\
    \le &  3\varepsilon \displaystyle\int_{\Omega} \eta^2 |\tau_hV_p(Du)|^2 dx \notag\\
&+ C_{\varepsilon}(L)|h|^{2\alpha}[D \psi]^2_{B^{\alpha}_{\frac{p+2\beta}{p+1+\beta-q}, \infty}(B_R)} \biggl( \displaystyle\int_{B_{t'}}(1+|D u|)^{p+2\beta}dx \biggr)^{\frac{2q-p-2}{p+2\beta}} \notag\\
&+\dfrac{C_{\varepsilon}(L)}{(t-s)^2} |h|^2 \biggl( \displaystyle\int_{B_t}|D (u-\psi)|^{\frac{p+2\beta}{p+1+\beta-q}}dx \biggr)^{\frac{2(p+1+\beta-q)}{p+2\beta}}  \biggl( \displaystyle\int_{B_{t'}}(1+|D u|)^{p+2\beta}dx \biggr)^{\frac{2q-p-2}{p+2\beta}} \notag\\
&+ C_{\varepsilon}|h|^{2 \alpha} \biggl( \displaystyle\int_{B_{R}}g^\frac{p+2\beta}{p+\beta-q} dx \biggr)^{\frac{2(p+\beta-q)}{p+2\beta}} \biggl(  \displaystyle\int_{B_t}(1+|Du|)^{p+2\beta}dx\biggr)^{\frac{2q-p}{p+2\beta}} \notag\\
&+|h|^{2\alpha}  \left( \int_{B_R} g^\frac{p+2\beta}{p+\beta-q} dx \right)^\frac{p+\beta-q}{p+2\beta}
[D \psi]_{B^\alpha_{\frac{p+2\beta}{p+1+\beta-q}, \infty}(B_R)}
\left( \ibt(1+|Du|)^{p+2\beta} dx \right)^\frac{q-1}{p+2\beta} \notag\\
&+\frac{C}{t-s} |h|^{\alpha+1}  \left( \int_{B_R} g^\frac{p+2\beta}{p+\beta-q} dx \right)^\frac{p+\beta-q}{p+2\beta}
\left( \ibtt |D (u- \psi)|^\frac{p+2\beta}{p+1+\beta-q} dx \right)^\frac{p+1+\beta-q}{p+2\beta} \notag\\
& \cdot 
\left( \ibt(1+|Du|)^{p+2\beta} dx \right)^\frac{q-1}{p+2\beta} .
\end{align}
Choosing $\varepsilon= \frac{C(\nu)}{6}$, we can reabsorb the first integral in the right hand side of the previous estimate
by the left hand side, thus getting

\begin{align}
      \displaystyle\int_{\Omega} &  \eta^{2} |\tau_{h}V_p(Du)|^{2} dx \notag\\
    \le &  
 C(L)|h|^{2\alpha}M_R^2 \biggl( \displaystyle\int_{B_{t'}}(1+|D u|)^{p+2\beta}dx \biggr)^{\frac{2q-p-2}{p+2\beta}} \notag\\
&+\dfrac{C(L)}{(t-s)^2} |h|^2 \biggl( \displaystyle\int_{B_t}|D (u-\psi)|^{\frac{p+2\beta}{p+1+\beta-q}}dx \biggr)^{\frac{2(p+1+\beta-q)}{p+2\beta}}  \biggl( \displaystyle\int_{B_{t'}}(1+|D u|)^{p+2\beta}dx \biggr)^{\frac{2q-p-2}{p+2\beta}} \notag\\
&+ C|h|^{2 \alpha} M_R^2 \biggl(  \displaystyle\int_{B_t}(1+|Du|)^{p+2\beta}dx\biggr)^{\frac{2q-p}{p+2\beta}} \notag\\
&+|h|^{2\alpha}2M_R
\left( \ibt(1+|Du|)^{p+2\beta} dx \right)^\frac{q-1}{p+2\beta} \notag\\
&+\frac{C}{t-s} |h|^{\alpha+1} M_R
\left( \ibtt |D (u- \psi)|^\frac{p+2\beta}{p+1+\beta-q} dx \right)^\frac{p+1+\beta-q}{p+2\beta}  
\left( \ibt(1+|Du|)^{p+2\beta} dx \right)^\frac{q-1}{p+2\beta} , \label{stimatauhvp}
\end{align}
where we set $M_R := \Vert g \Vert_{L^{\frac{p+2\beta}{p+\beta-q}}(B_R)}+ \Vert D\psi \Vert_{B^\alpha_{\frac{p+2\beta}{p+1+\beta-q},\infty}(B_R)}$.\\
From Young's inequality, we infer
\begin{align}
      \displaystyle\int_{\Omega} &  \eta^{2} |\tau_{h}V_p(Du)|^{2} dx \notag\\
    \le &  
 C_\theta(L,M_R)|h|^{2\alpha} + \theta |h|^{2\alpha} \displaystyle\int_{B_{t'}}(1+|D u|)^{p+2\beta}dx  \notag\\
&+\dfrac{C_\theta(L)}{(t-s)^{\tilde{p}}} |h|^2  \displaystyle\int_{B_t}|D (u-\psi)|^{\frac{p+2\beta}{p+1+\beta-q}}dx + \theta |h|^2  \displaystyle\int_{B_{t'}}(1+|D u|)^{p+2\beta}dx \notag\\
&+ C_\theta(M_R)|h|^{2 \alpha}  + \theta |h|^{2 \alpha}  \displaystyle\int_{B_t}(1+|Du|)^{p+2\beta}dx \notag\\
&+C_\theta(M_R) |h|^{2\alpha}  
+\theta |h|^{2\alpha}
 \ibt(1+|Du|)^{p+2\beta} dx  \notag\\
&+\frac{C_\theta(M_R)}{(t-s)^{p^*}} |h|^{\alpha+1}  
\left( \ibtt |D (u- \psi)|^\frac{p+2\beta}{p+1+\beta-q} dx \right)^\frac{p^*(p+1+\beta-q)}{p+2\beta} \notag\\
& + \theta |h|^{\alpha+1}  
\ibt(1+|Du|)^{p+2\beta} dx , \label{stimatau}
\end{align}
where $\tilde{p}:= \frac{p+2\beta}{p+\beta-q}$ and $p^*:= \frac{p+2\beta}{p+2\beta-q+1}$.\\
Using Young's inequality, we estimate the third and the second last integral appearing in the right hand side of estimate \eqref{stimatau} as follows 

\begin{align}
    & \dfrac{C_\theta(L)}{(t-s)^{\tilde{p}}}  |h|^2  \displaystyle\int_{B_t}|D (u-\psi)|^{\frac{p+2\beta}{p+1+\beta-q}}dx \notag\\
    & \le \dfrac{C_\theta(L)}{(t-s)^{\tilde{p}}} |h|^2  \displaystyle\int_{B_t}|D u|^{\frac{p+2\beta}{p+1+\beta-q}}dx
    + \dfrac{C_\theta(L)}{(t-s)^{\tilde{p}}} |h|^2  \displaystyle\int_{B_t}|D \psi|^{\frac{p+2\beta}{p+1+\beta-q}}dx \notag\\
    & \le \theta |h|^2  \displaystyle\int_{B_t}|D u|^{p+2\beta}dx + \dfrac{C_\theta(L)}{(t-s)^{p'}}  |h|^2 
    |B_R|
    + \dfrac{C_\theta(L,M_R)}{(t-s)^{\tilde{p}}} |h|^2 , \label{interpol1}
\end{align}
and analogously

\begin{align}
    &\frac{C_\theta(M_R)}{(t-s)^{p^*}} |h|^{\alpha+1} 
\left( \ibtt |D (u- \psi)|^\frac{p+2\beta}{p+1+\beta-q} dx \right)^\frac{p^*(p+1+\beta-q)}{p+2\beta} \notag\\
 & \le  C_\theta(M_R)  |h|^{\alpha+1} 
 + \frac{C_\theta}{(t-s)^{p''}} |h|^{\alpha+1}
 \ibtt |D (u- \psi)|^\frac{p+2\beta}{p+1+\beta-q} dx \notag\\
 & \le  C_\theta(M_R)  |h|^{\alpha+1}  
 + \theta |h|^2  \displaystyle\int_{B_t}|D u|^{p+2\beta}dx + \dfrac{C_\theta(L)}{(t-s)^{\tilde{p}}}  |h|^2 
    |B_R|  + \dfrac{C_\theta(L,M_R)}{(t-s)^{{p''}}} |h|^2  , \label{interpol2}
\end{align}
where $p':= \frac{p+1+\beta-q}{p+\beta-q}$, $p'':= \frac{p+2\beta}{p+1+\beta-q}$.\\
Inserting \eqref{interpol1} and \eqref{interpol2} in \eqref{stimatau}, we get

\begin{align}
      \displaystyle\int_{\Omega} &  \eta^{2} |\tau_{h}V_p(Du)|^{2} dx \notag\\
    \le &  
 C_\theta(L,M_R)|h|^{2\alpha} + \theta |h|^{2\alpha} \displaystyle\int_{B_{t'}}(1+|D u|)^{p+2\beta}dx  \notag\\
    & + \theta |h|^2  \displaystyle\int_{B_t}|D u|^{p+2\beta}dx + \dfrac{C_\theta(L)}{(t-s)^{p'}}  |h|^2 
    |B_R|
    + \dfrac{C_\theta(L,M_R)}{(t-s)^{\tilde{p}}} |h|^2  \notag\\
& + \theta |h|^2  \displaystyle\int_{B_{t'}}(1+|D u|)^{p+2\beta}dx \notag\\
&+ C_\theta(M_R)|h|^{2 \alpha}   + \theta |h|^{2 \alpha}  \displaystyle\int_{B_t}(1+|Du|)^{p+2\beta}dx \notag\\
&+C_\theta(M_R) |h|^{2\alpha}  
+\theta |h|^{2\alpha}
 \ibt(1+|Du|)^{p+2\beta} dx  \notag\\
& +C_\theta(M_R)  |h|^{\alpha+1} 
 + \theta |h|^2  \displaystyle\int_{B_t}|D u|^{p+2\beta}dx + \dfrac{C_\theta(L)}{(t-s)^{\tilde{p}}}  |h|^2 
    |B_R| \notag\\
    &  + \dfrac{C_\theta(L,M_R)}{(t-s)^{{p''}}} |h|^2 
 + \theta |h|^{\alpha+1}  
\ibt(1+|Du|)^{p+2\beta} dx . \
\end{align}
We can rewrite the previous estimate as

\begin{align}
    \displaystyle\int_{\Omega} &  \eta^{2} |\tau_{h}V_p(Du)|^{2} dx \notag\\
    \le &  5\theta|h|^{2 \alpha} \int_{B_t} (1+|Du|^{p+2\beta}) dx + 2\theta|h|^{2\alpha} \int_{B_{t'}} (1+|Du|^{p+2\beta}) dx \notag\\
    &+ C_\theta|h|^{2\alpha} \biggl( 1+
     \dfrac{1}{(t-s)^{p'}} 
    +\dfrac{1}{(t-s)^{\tilde{p}}}
    + \dfrac{1}{(t-s)^{p''}} \biggr),  \notag
\end{align}
for a constant $C_\theta:=C_\theta(\nu,L,M_R,R)$.\\
Dividing both sides of previous estimate by $|h|^{2\alpha}$, recalling that $\eta =1$ in $B_s$ and passing to the limit as $t' \rightarrow t^+$, we get

\begin{align}
    \displaystyle\int_{B_s} &  \dfrac{|\tau_{h}V_p(Du)|^{2}}{|h|^{2 \alpha}} dx \notag\\
    \le &  7 \theta\int_{B_{t}} (1+|Du|^{p+2\beta}) dx 
     \notag\\
    &+ C_\theta \biggl( 1 + \dfrac{1}{(t-s)^{p'}} 
    +\dfrac{1}{(t-s)^{\tilde{p}}}
    + \dfrac{1}{(t-s)^{p''}} \biggr) ,\label{3.21}
\end{align}
for every $h \in \R^n$.\\
Since $u  \in L^\infty_{loc}(\Omega)$ and $V_p(Du) \in B^\alpha_{2,\infty,loc}(\Omega)$, by virtue of Proposition \ref{mainint} and Theorem \ref{limiteq}, we infer the following inequality
\begin{align}
    \int_{B_\rho} |Du|^{p+2\beta}dx 
    \le &  C \Vert u \Vert^{2\beta}_{L^\infty(B_s)}
    \sup_h
    \displaystyle\int_{B_s}  \dfrac{|\tau_{h}V_p(Du)|^{2}}{|h|^{2 \alpha}} dx \notag\\
    &
     + 
    \frac{C}{(s-\rho)^{2p}}
    \Vert u \Vert^{2\beta}_{L^\infty(B_s)}
    \Vert u \Vert^{p}_{W^{1,p}(B_s)} \notag\\
    \le & 
    C (\Vert\psi\Vert_{L^\infty(B_R)}+\Vert u \Vert_{W^{1,p}(B_{R})})^\pi
    \sup_h
    \displaystyle\int_{B_s}  \dfrac{|\tau_{h}V_p(Du)|^{2}}{|h|^{2 \alpha}} dx \notag\\
    &+    \frac{C}{(s-\rho)^{2p}}  
    C (\Vert\psi\Vert_{L^\infty(B_R)}+\Vert u \Vert_{W^{1,p}(B_{R})})^\pi
    , \label{interpolazione}
\end{align}
for a constant $\pi :=\pi(n,p,q)$.\\
Taking the supremum over $h$ in the left hand side of \eqref{3.21} and using estimate \eqref{interpolazione}, we obtain

\begin{align}
     \int_{B_\rho} |Du|^{p+2\beta}dx 
    \le & 
    C (\Vert\psi\Vert_{L^\infty(B_R)}+\Vert u \Vert_{W^{1,p}(B_{R})})^\pi
     \theta\int_{B_t} (1+|Du|^{p+2\beta}) dx 
     \notag\\
    &+ C_\theta (\Vert\psi\Vert_{L^\infty(B_R)}+\Vert u \Vert_{W^{1,p}(B_{R})})^\pi
    \biggl[
    1 + \dfrac{1}{(t-s)^{p'}} 
    +\dfrac{1}{(t-s)^{\tilde{p}}}
    + \dfrac{1}{(t-s)^{p''}} \biggr]
    \notag\\
    &+    \frac{C}{(s-\rho)^{2p}}  
     (\Vert\psi\Vert_{L^\infty(B_R)}+\Vert u \Vert_{W^{1,p}(B_{R})})^\pi
    ,\label{3.23}
\end{align}
for every $0 < \frac{R}{4} < \rho  < s < t   <  \frac{R}{2} $.\\
Now, choosing $s$ such that $s-\rho=t-s$, i.e. $s=\frac{t+\rho}{2}$, it follows

\begin{align}
\int_{B_\rho} |Du|^{p+2\beta}dx 
    \le & 
    C (\Vert\psi\Vert_{L^\infty(B_R)}+\Vert u \Vert_{W^{1,p}(B_{R})})^\pi
     \theta\int_{B_t} (1+|Du|^{p+2\beta}) dx 
     \notag\\
    &+ C_\theta (\Vert\psi\Vert_{L^\infty(B_R)}+\Vert u \Vert_{W^{1,p}(B_{R})})^\pi
    \biggl[
    1 + \dfrac{1}{(t-\rho)^{p'}} 
    +\dfrac{1}{(t-\rho)^{\tilde{p}}}
    + \dfrac{1}{(t-\rho)^{p''}} \biggr]
    \notag\\
    &+    \frac{C}{(t-\rho)^{2p}}  
     (\Vert\psi\Vert_{L^\infty(B_R)}+\Vert u \Vert_{W^{1,p}(B_{R})})^\pi.
   \label{3.24}
\end{align}
Setting
\begin{gather}
\Phi (r) =  \displaystyle\int_{B_r}  |Du|^{p+2\beta} dx , \notag
\end{gather}
we can write inequality \eqref{3.24} as
\begin{align*}
\Phi(\rho)
    \le & 
    C (\Vert\psi\Vert_{L^\infty(B_R)}+\Vert u \Vert_{W^{1,p}(B_{R})})^\pi
     \theta\Phi(t)
     \notag\\
    &+ C_\theta (\Vert\psi\Vert_{L^\infty(B_R)}+\Vert u \Vert_{W^{1,p}(B_{R})})^\pi
    \biggl[
    1 + \dfrac{1}{(t-\rho)^{p'}} 
    +\dfrac{1}{(t-\rho)^{\tilde{p}}}
    + \dfrac{1}{(t-\rho)^{p''}} \biggr]
    \notag\\
    &+    \frac{C}{(t-\rho)^{2p}}  
     (\Vert\psi\Vert_{L^\infty(B_R)}+\Vert u \Vert_{W^{1,p}(B_{R})})^\pi.
\end{align*}
By virtue of Lemma \ref{lm2}, choosing $\theta $ such that $ C (\Vert\psi\Vert_{L^\infty(B_R)}+\Vert u \Vert_{W^{1,p}(B_{R})})^\pi
     \theta = \frac{1}{2}$, we obtain
\begin{align}
\Phi\biggl( \frac{R}{4} \biggr)
    \le & 
   C (\Vert\psi\Vert_{L^\infty(B_R)}+\Vert u \Vert_{W^{1,p}(B_{R})})^\pi
    \biggl[
    1 + \dfrac{1}{R^{p'}} 
    +\dfrac{1}{R^{\tilde{p}}}
    + \dfrac{1}{R^{p''}} \biggr]
    \notag\\
    &+    \frac{C}{R^{2p}}  
     (\Vert\psi\Vert_{L^\infty(B_R)}+\Vert u \Vert_{W^{1,p}(B_{R})})^\pi, \notag
\end{align}
with $C := C(n,p,q,\nu, L ,M_R)$.
\\Now, recalling the definition of $\Phi$, we obtain
\begin{align}
    \int_{B_{R/4}} |Du|^{p+2\beta} dx 
     \leq &
     C (\Vert\psi\Vert_{L^\infty(B_R)}+\Vert u \Vert_{W^{1,p}(B_{R})})^\pi \notag\\
     & \cdot
     \biggl(
     \int_{B_{R}} (k^{\frac{p+2\beta}{p+\beta-q}}+1) dx +\Vert D\psi \Vert_{B^\alpha_{\frac{p+2\beta}{p+1+\beta-q},\infty}(B_R)}
     \biggr)^\pi , \label{3.28}
\end{align}
thus, inserting \eqref{3.28} in \eqref{3.21}, we deduce the a priori estimate
\begin{align*}
\displaystyle\int_{B_{R/4}} |\tau_h V_p(Du)|^2 dx 
     \le &  C  |h|^{2 \alpha} (\Vert\psi\Vert_{L^\infty(B_R)}+\Vert u \Vert_{W^{1,p}(B_{R})})^\pi \notag\\
     & \cdot
     \biggl(
     \int_{B_{R}} (g^{\frac{p+2\beta}{p+\beta-q}}+1) dx +\Vert D\psi \Vert_{B^\alpha_{\frac{p+2\beta}{p+1+\beta-q},\infty}(B_R)}
     \biggr)^\pi ,
\end{align*}
for constants $C := C(n,p,q, \nu, L,R)$ and $\pi:= \pi (n,p,q,\beta)$. 
\endproof

Now, we are able to establish the following higher differentiability result for obstacle problems with $p$-growth.

\begin{thm}\label{thmp=q}
Assume that $\mathcal{A}(x,\xi)$ satisfies (A1)-(A3) for an exponent $2 \leq p=q $ and let $u \in \mathcal{K}_{\psi}(\Omega)$ be the solution to the obstacle problem \eqref{1}. If there exists a non-negative function $g \in L^{\infty}_{\text{loc}}(\Omega)$ such that
$$|\mathcal{A}(x,\xi)-\mathcal{A}(y, \xi)| \leq |x-y|^{\alpha} (g(x)+g(y)) (\mu^2 +|\xi|^2)^{\frac{p-1}{2}} ,$$
for a.e. $x,y \in \Omega$ and for every $\xi \in \mathbb{R}^n$, then the following implication
$$ \psi \in L^\infty_{loc}(\Omega), D\psi \in B^\alpha_{\frac{p+2\beta}{1+\beta},\infty,\textrm{loc}}(\Omega)\Rightarrow (\mu^2 +|Du|^2)^\frac{p-2}{4} Du \in B^\alpha_{2,\infty,\textrm{loc}}(\Omega),$$
holds, provided $0 < \beta < \alpha <1$.
\end{thm}
\proof Using Proposition \ref{lemma2.14}, we infer $D \psi \in L^{p+2\beta}_{loc}(\Omega)$. Hence, \cite[Theorem 2.6]{byun} yields
$Du \in L^{p+2\beta}_{loc}(\Omega).$
Arguing as in the proof of Theorem \ref{approximation}, we derive estimate \eqref{3.21} in the case $p=q$. This completes the proof.
\endproof

\subsection{Passage to the limit}\label{limit}
{\it Proof of Theorem \ref{mainthm1}.}
Let $u \in \mathcal{K}_{\psi}(\Omega)$ be a solution to \eqref{!1}, and let $F_j$ be defined as in Lemma \ref{apprlem1}. 
Fixed $B_R \Subset \Omega$, let $u_j$ be the solution of the problem
\begin{equation}
    \min \biggl\{ \displaystyle\int_{B_R} F_j(x,Dw) dx :w \geq \psi \ \text{a.e. in} \ B_R, \ w \in u + W^{1,p}_0(B_R)  \biggr\}. \label{10.2.}
\end{equation}

Setting 
\begin{gather}
\mathcal{A}_j(x,\xi)= D_{\xi}F_j(x,\xi), \notag
\end{gather}
one can easily check that $\mathcal{A}_j$ satisfies (A1)--(A4) and the following assumptions:
\begin{gather}
|\mathcal{A}_j(x,\xi)| \leq l_1(j)(\mu^2 +|\xi|^2)^{\frac{p-1}{2}} \label{K2} \\
|\mathcal{A}_j(x,\xi)-\mathcal{A}_j(x,\eta)|\leq L_1(j)|\xi-\eta| (\mu^2+|\xi|^{2}+|\eta|^{2})^{\frac{p-2}{2}} \label{K1}\\
|\mathcal{A}_j(x,\xi)-\mathcal{A}_j(y, \xi)| \leq \Theta(j)|x-y|^{\alpha} (g(x)+g(y)) (\mu^2 +|\xi|^2)^{\frac{p-1}{2}} \label{K3}
\end{gather}
for a.e. $x,y \in \Omega$, for every $\xi, \eta \in \mathbb{R}^n$ and $j \in \N$. It is well known that $u_j \in \mathcal{K}_{\psi}(B_R)$ is a minimizer of problem \eqref{10.2.} if, and only if, the following variational inequality holds 
\begin{align}\label{varineq2.1}
\displaystyle\int_{\Omega} \langle \mathcal{A}_j(x,Du_j)  ,D(\varphi-u_j) \rangle dx \geq 0 , \quad \forall \varphi\in \mathcal{K}_{\psi}(B_R).
\end{align}

Let $\Omega' \Subset \Omega$ be an open set. Fix a non-negative smooth kernel $\phi \in \mathcal{C}^\infty_0(B_1(0))$ such that $\int_{B_1(0)} \phi =1$ and consider the corresponding family of mollifiers $\{\phi_m\}_{m \in \N}$. Setting 
$$g_m=g * \phi_m$$
and
\begin{gather}
\mathcal{A}_{jm}(x,\xi)= \int_{B_1(0)} \phi(y)  \mathcal{A}_j(x+my, \xi) dy, \label{A_j}
\end{gather}
an easy computation shows that $\mathcal{A}_{jm}$ satisfies assumptions (A1)--(A3), \eqref{K2}--\eqref{K1} and the conditions:
$$|\mathcal{A}_{jm}(x,\xi)-\mathcal{A}_{jm}(y, \xi)| \leq |x-y|^{\alpha} (g_m(x)+g_m(y)) (\mu^2 +|\xi|^2)^{\frac{q-1}{2}} \eqno{\rm{{ (A4')}}}$$
\begin{equation} \label{K4}
|\mathcal{A}_{jm}(x,\xi)-\mathcal{A}_{jm}(y, \xi)| \leq \Theta(j)|x-y|^{\alpha} (g_m(x)+g_m(y)) (\mu^2 +|\xi|^2)^{\frac{p-1}{2}} 
\end{equation}
for a.e. $x,y \in \Omega$, for every $\xi, \eta \in \mathbb{R}^n$ and every $j ,m \in \N$.

\textbf{Step 1.} Fixed $j \in \N$,
let $\mathcal{A}_{jm}$ be defined as in \eqref{A_j} and let $u_{jm} \in u_j+W_0^{1,p}(B_R)$ be the solution to the variational inequality 

\begin{equation}\label{10.1}
    \int_{B_R} \langle \mathcal{A}_{jm}(x,Du_{jm}), D(\varphi-u_{jm}) \rangle dx \ge 0, \quad \forall \varphi \in \mathcal{K}_{\psi}(B_R).
\end{equation}
By the ellipticity assumption (A1), we have
\begin{align}
    \nu & \int_{B_R}  (\mu^2+|Du_j|^2+|Du_{jm}|^2)^{\frac{p-2}{2}}
    |Du_{jm}-Du_j|^2 dx \notag\\
    \leq & \int_{B_R}
    \langle \mathcal{A}_{jm}(x,Du_{jm})
    -\mathcal{A}_{jm}(x,Du_j),
    Du_{jm}-Du_j\rangle dx \notag\\
    =& \int_{B_R}
    \langle \mathcal{A}_{jm}(x,Du_{jm}),
    Du_{jm}-Du_j\rangle dx \notag\\
    &- \int_{B_R}
    \langle
    \mathcal{A}_{jm}(x,Du_j),
    Du_{jm}-Du_j\rangle dx \notag\\
    =& \int_{B_R}
    \langle \mathcal{A}_{jm}(x,Du_{jm}),
    Du_{jm}-Du_j\rangle dx \notag\\
    &- \int_{B_R}
    \langle
    \mathcal{A}_j(x,Du_j),
    Du_{jm}-Du_j\rangle dx \notag\\
    &+ \int_{B_R}
    \langle \mathcal{A}_j(x,Du_{j})
    -\mathcal{A}_{jm}(x,Du_j),
    Du_{jm}-Du_j\rangle dx \label{3.2.1}
\end{align}
Since $u_{j}$ and $u_{jm}$ are solutions to \eqref{10.2.} and \eqref{10.1} respectively, we notice that

\begin{align}
   & \int_{B_R}
    \langle \mathcal{A}_{jm}(x,Du_{jm}),
    Du_{jm}-Du_j\rangle dx 
    - \int_{B_R}
    \langle
    \mathcal{A}_j(x,Du_j),
    Du_{jm}-Du_j\rangle dx \le 0 \label{3.2.2}
\end{align}
Combining \eqref{3.2.1} and \eqref{3.2.2}, we get

\begin{align}
    \nu & \int_{B_R} (\mu^2+|Du_j|^2+|Du_{jm}|^2)^{\frac{p-2}{2}}
    |Du_{jm}-Du_j|^2 dx \notag\\
    \le & \int_{B_R}
    \langle \mathcal{A}_j(x,Du_{j})
    -\mathcal{A}_{jm}(x,Du_j),
    Du_{jm}-Du_j\rangle dx \notag\\
    \le & \biggl( \int_{B_R}
    | \mathcal{A}_j(x,Du_{j})
    -\mathcal{A}_{jm}(x,Du_j)|^\frac{p}{p-1} dx  \biggr)^\frac{p-1}{p}
    \biggl( \int_{B_R}|
    Du_{jm}-Du_j|^p dx  \biggr)^\frac{1}{p}, \label{3.2.3}
\end{align}
where in the last inequality we used H\"{o}lder's inequality.\\
Since $p \ge 2$, from \eqref{3.2.3}, we obtain
\begin{align} \label{3.2.5}
    \int_{B_R}|
    Du_{jm}-Du_j|^p dx
    \le C \int_{B_R}
    | \mathcal{A}(x,Du_{j})
    -\mathcal{A}_{jm}(x,Du_j)|^\frac{p}{p-1} dx.
\end{align}
Since $\mathcal{A}_{jm}(x,Du_{j})$ satisfies 
$$|\mathcal{A}_{jm}(x,Du_j)| \leq l_1(j)(\mu^2 +|Du_j|^2)^{\frac{p-1}{2}},$$
and $\mathcal{A}_{jm}(x,Du_{j}) \rightarrow_{m \rightarrow \infty} \mathcal{A}(x,Du_j)$ a.e. in $\Omega$, applying the Dominated convergence theorem, we have
$$  \mathcal{A}_{jm}(x,Du_{j}) \rightarrow_{m \rightarrow \infty }\mathcal{A}(x,Du_j)  \quad \text{strongly in } L^\frac{p}{p-1}(\Omega) .$$
Therefore, passing to the limit for $m \rightarrow\infty$ in \eqref{3.2.5}, we deduce that 
\begin{equation}
    u_{jm} \rightarrow u_j \ \text{in} \ W^{1,p}(B_R). \label{strongconvergence}
\end{equation}
Moreover, since $g \in L^{\frac{p+2\beta}{p+\beta-q}}_{loc}(\Omega)$, we have
\begin{equation}
    g_{m} \rightarrow g \ \text{in} \ L^{\frac{p+2\beta}{p+\beta-q}}_{loc}(\Omega). \label{strongcon1}
\end{equation}
By virtue of Theorem \ref{thmp=q}, $V_p(Du_{jm}) \in B^\alpha_{2,\infty,loc}(B_R)$. Hence, from Theorem \ref{approximation}, $u_{jm} $ satisfies the a priori estimate 
\begin{align}
    \int_{B_{R/4}} |Du_{jm}|^{p+2\beta} dx 
     \leq &
     C (\Vert\psi\Vert_{L^\infty(B_R)}+\Vert u_{jm} \Vert_{W^{1,p}(B_{R})})^\pi \notag\\
     & \cdot
     \biggl(
     \int_{B_{R}} (g_m^{\frac{p+2\beta}{p+\beta-q}}+1) dx +\Vert D\psi \Vert_{B^\alpha_{\frac{p+2\beta}{p+1+\beta-q},\infty}(B_R)}
     \biggr)^\pi ,
     \end{align}
for constants $C := C(n,p,q, \nu, L,R)$ and $\pi:= \pi (n,p,q,\beta)$, both independent of $j$ and $m$.
\\Finally, by weak lower semicontinuity, \eqref{strongconvergence} and \eqref{strongcon1}, we get

\begin{align}
    \int_{B_{R/4}} |Du_{j}|^{p+2\beta} dx \le & \liminf_{m \rightarrow \infty} \int_{B_{R/4}} |Du_{jm}|^{p+2\beta} dx \notag\\
    \leq &
    C
     (\Vert\psi\Vert_{L^\infty(B_R)}+\Vert u_{j} \Vert_{W^{1,p}(B_{R})})^\pi \notag\\
     & \cdot  \biggl(
     \int_{B_{R}} (g^{\frac{p+2\beta}{p+\beta-q}}+1) dx +\Vert D\psi \Vert_{B^\alpha_{\frac{p+2\beta}{p+1+\beta-q},\infty}(B_R)}
     \biggr)^\pi. \label{10.3.}
\end{align}
\vspace{0.5cm}
\textbf{Step 2.} 
From Lemma \ref{apprlem1} $(iv)$, there exists $c_1 > 0$ such that
\begin{equation}
|\xi|^p \leq c_1 (1+ F_j(x,\xi)), \quad \forall j \in \mathbb{N}.\notag
\end{equation}
The previous estimate and the minimality of $u_j$ imply
\begin{align}\label{intF_j}
\displaystyle\int_{B_R}|Du_j|^p dx \leq & c_1 \displaystyle\int_{B_R} (1+ F_j(x,Du_j))dx \notag \\
\leq & c_1 \displaystyle\int_{B_R} (1+ F_j(x,Du))dx \notag \\
\leq & c_1 \displaystyle\int_{B_R} (1+ F(x,Du))dx, 
\end{align}
where in the last inequality we used Lemma \ref{apprlem1} $(ii)$. Thus, up to subsequences, 
\begin{equation}\label{weakconv}
u_j \rightharpoonup \tilde{u} \ \text{in} \ u+W^{1,p}_{0}(B_R).
\end{equation}
Now, fix $j_0 \in \mathbb{N}$. Then, by Lemma \ref{apprlem1} $(ii)$ and the fact that $u_j$ is a minimum for $F_j$, for every $j > j_0$, we might write
\begin{align*}
\int_{B_R} F_{j_0}(x,Du_j) dx &\leq \int_{B_R} F_j(x,Du_j)dx\\
&\leq \int_{B_R} F_j(x,Du)dx\leq \int_{B_R} F(x,Du)dx.
\end{align*}
From weak lower semicontinuity of $F_{j_0}$ and \eqref{weakconv}, it holds, 
\begin{align*}
\int_{B_R} F_{j_0}(x,D\tilde{u})dx\leq \liminf_{j\to +\infty}\int_{B_R} F_{j_0}(x,Du_j)dx \leq \int_{B_R} F(x,Du)dx .
\end{align*}
Combining these last inequalities, we get
\begin{align*}
\int_{B_R} F(x,D\tilde{u})dx= \lim_{j_0 \to +\infty}\int_{B_R} F_{j_0}(x,D\tilde{u})dx\leq \int_{B_R}F(x,Du)dx,
\end{align*}
where we also applied the monotone convergence theorem, according to Lemma \ref{apprlem1} $(ii)$.
\\Moreover, by the weak convergence \eqref{weakconv}, the limit function $\tilde{u}$ still belongs to $\mathcal{K}_{\psi}(B_R)$, since this set is convex and closed. Thus, by strict convexity of $F$, we have that $\tilde{u}= u $ a.e. in $B_R$. 

Now, from estimates \eqref{10.3.} and \eqref{intF_j}, it follows
\begin{align}
    \int_{B_{R/4}} |Du_j|^{p+2\beta} dx 
     \leq &
     C (\Vert\psi\Vert_{L^\infty(B_R)}+\Vert u_j \Vert_{W^{1,p}(B_{R})})^\pi \notag\\
     & \cdot
     \biggl(
     \int_{B_{R}} (g^{\frac{p+2\beta}{p+\beta-q}}+1) dx +\Vert D\psi \Vert_{B^\alpha_{\frac{p+2\beta}{p+1+\beta-q},\infty}(B_R)}
     \biggr)^\pi \notag\\
     \le & C \biggl(\Vert\psi\Vert_{L^\infty(B_R)}+\int_{B_R}(1+|u|^p+F(x,Du))dx \biggr)^\pi \notag\\
     & \cdot
     \biggl(
     \int_{B_{R}} (g^{\frac{p+2\beta}{p+\beta-q}}+1) dx +\Vert D\psi \Vert_{B^\alpha_{\frac{p+2\beta}{p+1+\beta-q},\infty}(B_R)}
     \biggr)^\pi, \notag
\end{align}
for constants $C := C(n,p,q, \nu, L,R)$ and $\pi:= \pi (n,p,q,\beta)$, both independent of $j$.
\\Hence, from  \eqref{weakconv} and weak lower semicontinuity, it follows
\begin{align}
    \int_{B_{R/4}} |Du|^{p+2\beta} dx \le & \liminf_{j \rightarrow \infty} \int_{B_{R/4}} |Du_{j}|^{p+2\beta} dx \notag\\
    \le & C \biggl(\Vert\psi\Vert_{L^\infty(B_R)}+\int_{B_R}(1+|u|^p+F(x,Du))dx \biggr)^\pi \notag\\
     & \cdot
     \biggl(
     \int_{B_{R}} (g^{\frac{p+2\beta}{p+\beta-q}}+1) dx +\Vert D\psi \Vert_{B^\alpha_{\frac{p+2\beta}{p+1+\beta-q},\infty}(B_R)}
     \biggr)^\pi. \notag
\end{align}
Eventually, proceeding as in the proof of Theorem \ref{approximation}, we derive that $V_p(Du) \in B^\alpha_{2,\infty,loc}(\Omega) $. 
\qed

\section{Proof of Theorem \ref{mainthm2}}\label{teorema1.2}
{\it Proof of Theorem \ref{mainthm2}.}
We derive only the a priori estimates, since the approximation procedure is achieved using the same arguments presented in Section \ref{limit}.\\
We a priori assume that $V_p(Du) \in B^\alpha_{2,\sigma,loc}(\Omega)$.
By virtue of assumption \eqref{gap2} and Theorem \ref{thmbound}, $u \in L^\infty_{loc}(\Omega)$.
Hence, using Proposition \ref{higherintfin}, we deduce that
$$Du \in L^{p+2 \alpha}_{loc}(\Omega).$$
Arguing analogously as in the proof of Theorem \ref{approximation}, but taking into account the new assumptions (A5) on the coefficients of the map $\mathcal{A}(x,\xi)$ and on the gradient of the obstacle, from inequality \eqref{stimatauhvp}, we obtain the following estimate

\begin{align}
      \displaystyle\int_{\Omega} &  \eta^{2} |\tau_{h}V_p(Du)|^{2} dx \notag\\
    \le &  
 C\biggl( \int_{B_R}|D\psi|^{\frac{p+2\alpha}{p+1+\alpha-q} } dx\biggr)^\frac{2 (p+1+\alpha -q)}{p+2\alpha}\biggl( \displaystyle\int_{B_{t'}}(1+|D u|)^{p+2\alpha}dx \biggr)^{\frac{2q-p-2}{p+2\alpha}} \notag\\
&+\dfrac{C}{(t-s)^2} |h|^2 \biggl( \displaystyle\int_{B_t}|D (u-\psi)|^{\frac{p+2\alpha}{p+1+\alpha-q}}dx \biggr)^{\frac{2(p+1+\alpha-q)}{p+2\alpha}}  \biggl( \displaystyle\int_{B_{t'}}(1+|D u|)^{p+2\alpha}dx \biggr)^{\frac{2q-p-2}{p+2\alpha}} \notag\\
&+ C|h|^{2 \alpha} \biggl( \displaystyle\int_{B_{R/2}}(g_k(x+h)+g_k(x))^\frac{p+2\alpha}{p+\alpha-q} dx \biggr)^{\frac{2(p+\alpha-q)}{p+2\alpha}} \biggl(  \displaystyle\int_{B_t}(1+|Du|)^{p+2\alpha}dx\biggr)^{\frac{2q-p}{p+2\alpha}} \notag\\
&+|h|^{\alpha}  \left( \int_{B_{R/2}} (g_k(x+h)+g_k(x))^\frac{p+2\alpha}{p+\alpha-q} dx \right)^\frac{p+\alpha-q}{p+2\alpha} \notag\\
& \cdot
\biggl( \int_{B_R}|\tau_hD \psi|^{\frac{p+2\alpha}{p+1+\alpha-q}} dx \biggr)^\frac{p+1+\alpha-q}{p+2\alpha}
\left( \ibt(1+|Du|)^{p+2\alpha} dx \right)^\frac{q-1}{p+2\alpha} \notag\\
&+\frac{C}{t-s} |h|^{\alpha+1}  \left( \int_{B_{R/2}} (g_k(x+h)+g_k(x))^\frac{p+2\alpha}{p+\alpha-q} dx \right)^\frac{p+\alpha-q}{p+2\alpha}
\left( \ibtt |D (u- \psi)|^\frac{p+2\alpha}{p+1+\alpha-q} dx \right)^\frac{p+1+\alpha-q}{p+2\alpha} \notag\\
& \cdot 
\left( \ibt(1+|Du|)^{p+2\alpha} dx \right)^\frac{q-1}{p+2\alpha} ,
\end{align}
for a positive constant $C:=C(n,p,q,\nu,L)$, where $2^{-k}(t'-t) \le |h| \le 2^{-k+1}(t'-t)$, $k \in \N$. Recalling that $\eta=1 $ on $B_s$ and dividing both sides by $|h|^{2 \alpha}$, we get

\begin{align}
      \displaystyle\int_{B_s} & \frac{|\tau_{h}V_p(Du)|^{2}}{|h|^{2\alpha}} dx \notag\\
    \le &  
 C\biggl( \int_{B_R}\frac{|D\psi|^{\frac{p+2\alpha}{p+1+\alpha-q} }}{|h|^\frac{\alpha(p+2\alpha)}{p+1+\alpha-q}} dx\biggr)^\frac{2 (p+1+\alpha -q)}{p+2\alpha}\biggl( \displaystyle\int_{B_{t'}}(1+|D u|)^{p+2\alpha}dx \biggr)^{\frac{2q-p-2}{p+2\alpha}} \notag\\
&+\dfrac{C}{(t-s)^2} |h|^{2(1-\alpha)} \biggl( \displaystyle\int_{B_t}|D (u-\psi)|^{\frac{p+2\alpha}{p+1+\alpha-q}}dx \biggr)^{\frac{2(p+1+\alpha-q)}{p+2\alpha}}  \biggl( \displaystyle\int_{B_{t'}}(1+|D u|)^{p+2\alpha}dx \biggr)^{\frac{2q-p-2}{p+2\alpha}} \notag\\
&+ C \biggl( \displaystyle\int_{B_{R/2}}(g_k(x+h)+g_k(x))^\frac{p+2\alpha}{p+\alpha-q} dx \biggr)^{\frac{2(p+\alpha-q)}{p+2\alpha}} \biggl(  \displaystyle\int_{B_t}(1+|Du|)^{p+2\alpha}dx\biggr)^{\frac{2q-p}{p+2\alpha}} \notag\\
&+  \left( \int_{B_{R/2}} (g_k(x+h)+g_k(x))^\frac{p+2\alpha}{p+\alpha-q} dx \right)^\frac{p+\alpha-q}{p+2\alpha} \notag\\
& \cdot
\biggl( \int_{B_R}\frac{|\tau_h D \psi|^{\frac{p+2\alpha}{p+1+\alpha-q}}}{|h|^\frac{\alpha(p+2\alpha)}{p+1+\alpha-q}} dx \biggr)^\frac{p+1+\alpha-q}{p+2\alpha}
\left( \ibt(1+|Du|)^{p+2\alpha} dx \right)^\frac{q-1}{p+2\alpha} \notag\\
&+\frac{C}{t-s} |h|^{1-\alpha}  \left( \int_{B_{R/2}} (g_k(x+h)+g_k(x))^\frac{p+2\alpha}{p+\alpha-q} dx \right)^\frac{p+\alpha-q}{p+2\alpha}
\left( \ibtt |D (u- \psi)|^\frac{p+2\alpha}{p+1+\alpha-q} dx \right)^\frac{p+1+\alpha-q}{p+2\alpha} \notag\\
& \cdot 
\left( \ibt(1+|Du|)^{p+2\alpha} dx \right)^\frac{q-1}{p+2\alpha} .
\end{align}
We need now to take the $L^\sigma$ norm with the measure $\frac{d h}{|h|^n}$ restricted to the ball $B(0,t'-t)$ on the $h$-space of the $L^2$ norm of the difference quotient of order $\alpha$ of the function $V_p(Du)$. Since the functions $g_k$ are defined for $2^{-k}(t'-t) \leq |h| \leq 2^{-k+1}(t'-t)$ we interpret the ball $B(0,t'-t)$ as
$$ B(0,t'-t)= \bigcup_{k=1}^{\infty} B(0,2^{-k+1}(t'-t))\setminus B(0,2^{-k}(t'-t))=: \bigcup_{k=1}^{\infty} E_k.$$
We obtain the following estimate

\begin{align}
& \int_{B_{t'-t}(0)} \biggl(
      \displaystyle\int_{B_s}  \frac{|\tau_{h}V_p(Du)|^{2}}{|h|^{2\alpha}} dx
      \biggr)^\frac{\sigma}{2} \frac{dh}{|h|^n}\notag\\
    \le &  
 C
  \biggl( \displaystyle\int_{B_{t'}}(1+|D u|)^{p+2\alpha}dx \biggr)^{\frac{\sigma(2q-p-2)}{2(p+2\alpha)}}
   \int_{B_{t'-t}(0)} \biggl(
  \int_{B_R}\frac{|D\psi|^{\frac{p+2\alpha}{p+1+\alpha-q} }}{|h|^\frac{\alpha(p+2\alpha)}{p+1+\alpha-q}} dx  \biggr)^\frac{\sigma (p+1+\alpha -q)}{p+2\alpha}
  \frac{dh}{|h|^n} 
  \notag\\
&+\dfrac{C}{(t-s)^2}\biggl( \displaystyle\int_{B_t}|D (u-\psi)|^{\frac{p+2\alpha}{p+1+\alpha-q}}dx \biggr)^{\frac{\sigma(p+1+\alpha-q)}{p+2\alpha}}  \biggl( \displaystyle\int_{B_{t'}}(1+|D u|)^{p+2\alpha}dx \biggr)^{\frac{(2q-p-2)\sigma}{2(p+2\alpha)}}
\int_{B_{t'-t}(0)}  |h|^{\sigma(1-\alpha)}   \frac{dh}{|h|^n} 
\notag\\
&+C  \biggl(  \displaystyle\int_{B_t}(1+|Du|)^{p+2\alpha}dx\biggr)^{\frac{(2q-p)\sigma}{2(p+2\alpha)}}
\sum_{k=1}^\infty\int_{E_k}  
\biggl( \displaystyle\int_{B_{R/2}}(g_k(x+h)+g_k(x))^\frac{p+2\alpha}{p+\alpha-q} dx \biggr)^{\frac{\sigma(p+\alpha-q)}{p+2\alpha}}
\frac{dh}{|h|^n}
\notag\\
&+ \sum_{k=1}^\infty\int_{E_k}  \left( \int_{B_{R/2}} (g_k(x+h)+g_k(x))^\frac{p+2\alpha}{p+\alpha-q} dx \right)^\frac{(p+\alpha-q)\sigma}{2(p+2\alpha)}
\biggl( \int_{B_R}\frac{|\tau_h D \psi|^\frac{p+2\alpha}{p+1+\alpha-q}}{|h|^\frac{\alpha(p+2\alpha)}{p+1+\alpha-q}} dx \biggr)^\frac{(p+1+\alpha-q)\sigma}{2(p+2\alpha)}
\frac{dh}{|h|^n} \notag\\
& \cdot
\left( \ibt(1+|Du|)^{p+2\alpha} dx \right)^\frac{(q-1)\sigma}{2(p+2\alpha)} \notag\\
&+\frac{C}{t-s} \sum_{k=1}^\infty\int_{E_k}|h|^{\frac{\sigma(1-\alpha)}{2}}  \left( \int_{B_{R/2}} (g_k(x+h)+g_k(x))^\frac{p+2\alpha}{p+\alpha-q} dx \right)^\frac{(p+\alpha-q)\sigma}{2(p+2\alpha)} \notag\\
& \cdot
\left( \ibtt |D (u- \psi)|^\frac{p+2\alpha}{p+1+\alpha-q} dx \right)^\frac{(p+1+\alpha-q)\sigma}{2(p+2\alpha)} 
\left( \ibt(1+|Du|)^{p+2\alpha} dx \right)^\frac{(q-1)\sigma}{2(p+2\alpha)}. \label{stimateorema2}
\end{align}
Note that, since $\alpha \leq \gamma$, the integral 
$$J_1:= \int_{B_{t'-t}(0)} \biggl( \int_{B_{R}}\dfrac{|\tau_h D \psi|^\frac{p+2\alpha}{p+1+\alpha-q}}{|h|^\frac{\alpha(p+2\alpha)}{p+1+\alpha-q}}dx \biggr)^{\frac{\sigma(p+1+\alpha-q)}{p+2\alpha}} \dfrac{d h}{|h|^n}  $$
is controlled by the norm in the Besov space $B^{\alpha}_{\frac{p+2\alpha}{p+1+\alpha-q},\sigma}$ on $B_{R}$ of the gradient of the obstacle which is finite by assumptions. The integral
$$ J_2 := \int_{B_{t'-t}(0)} |h|^{\sigma(1-\alpha)} \dfrac{dh}{|h|^n}  $$ 
can be calculated in polar coordinates as follows
$$J_2 = C(n) \displaystyle\int_0^{t'-t} \varrho^{\sigma(1-\alpha)-1} d \varrho\le C(n) \displaystyle\int_0^{R/4} \varrho^{\sigma(1-\alpha)-1} d \varrho = C(n,\alpha,\sigma,R),$$
since $t'-t \le \frac{R}{4}$ and $\alpha \in (0,1)$.
\\Now, we take care of the integral
$$J_3:= \displaystyle\sum_{k=1}^{\infty}
\int_{E_k}
 \biggl( \displaystyle\int_{B_{R/2}}(g_k(x+h)+g_k(x))^\frac{p+2\alpha}{p+\alpha-q} dx \biggr)^{\frac{\sigma(p+\alpha-q)}{p+2\alpha}} \dfrac{dh}{|h|^n} .$$
We write the right hand sinde of the previous estimate in polar coordinates, so $h \in E_k$ if, and only if, $h= r \xi$ for some $2^{-k+1}(t'-t) \leq m < 2^{-k}(t'-t)$ and some $\xi$ in the unit sphere $\mathbb{S}^{n-1}$ on $\R^n$. We denote by $d S(\xi)$ the surface measure on $\mathbb{S}^{n-1}$. We infer
\begin{align*}
J_3 \leq &  \displaystyle\sum_{k=1}^{\infty} \displaystyle\int_{m_{k-1}}^{m_k} \displaystyle\int_{\mathbb{S}^{n-1}}  \biggl( \displaystyle\int_{B_{R/2}}(g_k(x+h)+g_k(x))^\frac{p+2\alpha}{p+\alpha-q} dx \biggr)^{\frac{\sigma(p+\alpha-q)}{p+2\alpha}} dS(\xi)\dfrac{dm}{m} \notag\\
= & \displaystyle\sum_{k=1}^{\infty} \displaystyle\int_{m_{k-1}}^{m_k} \displaystyle\int_{\mathbb{S}^{n-1}}  \Vert (\tau_{m \xi}g_k+g_k) \Vert_{L^\frac{p+2\alpha}{p+\alpha-q}(B_{R/2})}^{\sigma} dS(\xi)\dfrac{dm}{m},
\end{align*}
where we set $m_k=2^{-k}(t'-t)$.  We note that for each $\xi \in \mathbb{S}^{n-1}$ and $m_{k-1}\leq m \leq m_k$
\begin{align*}
\Vert (\tau_{m \xi}g_k+g_k) \Vert_{L^\frac{p+2\alpha}{p+\alpha-q}(B_{R/2})}  \leq & \Vert g_k \Vert_{L^\frac{p+2\alpha}{p+\alpha-q}(B_{R/2}-m_k \xi)} + \Vert g_k \Vert_{L^\frac{p+2\alpha}{p+\alpha-q}(B_{R/2})} \notag\\
\leq & 2  \Vert g_k \Vert_{L^\frac{p+2\alpha}{p+\alpha-q}(B_{R/2+R/4})},
\end{align*}
where in the last inequality we used that $t'-t \le \frac{R}{4}$.
Hence
$$J_3 \leq C(n) \Vert \{g_k \}_k \Vert_{l^\sigma(L^\frac{p+2\alpha}{p+\alpha-q}(B_{R}))}^{\sigma},$$
which is finite by assumption (F6).
\\Using the Young's inequality with exponent $2$, we deduce the following estimate
\begin{align*}
& \displaystyle\sum_{k=1}^{\infty}\int_{E_k}  \left( \displaystyle\int_{B_{R/2}} (g_k(x+h)+g_k(x))^\frac{p+2\alpha}{p+\alpha-q}  dx \right)^\frac{\sigma(p+\alpha-q)}{2(p+2\alpha)} \left( \displaystyle\int_{B_{R/2}} \dfrac{|\tau_h D\psi|^\frac{p+2\alpha}{p+1+\alpha-q}}{|h|^{ \frac{\alpha(p+2\alpha)}{p+1+\alpha-q}}} dx \right)^\frac{\sigma(p+1+\alpha-q)}{2(p+2\alpha)} \dfrac{dh}{|h|^n} \\
\leq & C \displaystyle\sum_{k=1}^{\infty}\int_{E_k}  \left( \displaystyle\int_{B_{R/2}} (g_k(x+h)+g_k(x))^\frac{p+2\alpha}{p+\alpha-q}  dx \right)^\frac{\sigma(p+\alpha-q)}{p+2\alpha} \dfrac{dh}{|h|^n} \notag\\
&+ C\int_{B_{t'-t}(0)} \left( \displaystyle\int_{B_{R/2}} \dfrac{|\tau_h D\psi|^\frac{p+2\alpha}{p+1+\alpha-q}}{|h|^{ \frac{\alpha(p+2\alpha)}{p+1+\alpha-q}}} dx \right)^\frac{\sigma(p+1+\alpha-q)}{p+2\alpha} \dfrac{dh}{|h|^n} 
\end{align*} 
where the two integrals in the right hand side can be estimated as the integrals $J_1$ and $J_3$.
\\Similarly, we obtain
\begin{align*}
& \displaystyle\sum_{k=1}^{\infty}\int_{E_k} |h|^{(1- \alpha)\frac{\sigma}{2}} \left(\displaystyle\int_{B_{R/2}} (g_k(x+h)+g_k(x))^\frac{p+2\alpha}{p+\alpha-q} dx \right)^\frac{\sigma(p+\alpha-q)}{2(p+2\alpha)} \dfrac{dh}{|h|^n}\\
\leq &\int_{B_{t'-t}(0)} |h|^{(1- \alpha)\sigma} \dfrac{dh}{|h|^n} + \displaystyle\sum_{k=1}^{\infty}\int_{E_k} \left(\displaystyle\int_{B_{R/2}} (g_k(x+h)+g_k(x))^\frac{p+2\alpha}{p+\alpha-q} dx \right)^\frac{\sigma(p\alpha-q)}{p+2\alpha} \dfrac{dh}{|h|^n}.
\end{align*}
The first term and the latter one can be estimated as the integral $J_2$ and $J_3$, respectively. 

Estimate \eqref{stimateorema2} can be written in the following way

\begin{align}
 \int_{B_{t'-t}(0)} & \biggl(
      \displaystyle\int_{B_s}  \frac{|\tau_{h}V_p(Du)|^{2}}{|h|^{2\alpha}} dx
      \biggr)^\frac{\sigma}{2} \frac{dh}{|h|^n}\notag\\
    \le &  
 C
  \biggl( \displaystyle\int_{B_{t'}}(1+|D u|)^{p+2\alpha}dx \biggr)^{\frac{\sigma(2q-p-2)}{2(p+2\alpha)}} \notag\\
  &
  +\dfrac{C}{(t-s)^2}
  \biggl( \displaystyle\int_{B_t}|D (u-\psi)|^{\frac{p+2\alpha}{p+1+\alpha-q}}dx \biggr)^{\frac{\sigma(p+1+\alpha-q)}{p+2\alpha}}
  \biggl( \displaystyle\int_{B_{t'}}(1+|D u|)^{p+2\alpha}dx \biggr)^{\frac{(2q-p-2)\sigma}{2(p+2\alpha)}}
\notag\\
&+ C  \biggl(  \displaystyle\int_{B_t}(1+|Du|)^{p+2\alpha}dx\biggr)^{\frac{(2q-p)\sigma}{2(p+2\alpha)}}
+ C
\left( \ibt(1+|Du|)^{p+2\alpha} dx \right)^\frac{(q-1)\sigma}{2(p+2\alpha)} \notag\\
&+\frac{C}{t-s} 
\left( \ibtt |D (u- \psi)|^\frac{p+2\alpha}{p+1+\alpha-q} dx \right)^\frac{(p+1+\alpha-q)\sigma}{2(p+2\alpha)} 
\left( \ibt(1+|Du|)^{p+2\alpha} dx \right)^\frac{(q-1)\sigma}{2(p+2\alpha)}, 
\end{align}
for a constant $C:=	C(L, \nu,p,q,r,n,\sigma,\alpha,R,\Vert D\psi \Vert_{B^{\gamma}_{\frac{p+2\alpha}{p+1+\alpha-q},\sigma}(B_{R})}, \Vert \{g_k \}_k \Vert_{l^\sigma(L^{\frac{p+2\alpha}{p+\alpha-q}}(B_{R}))})$.\\
From Young's inequality, we infer

\begin{align}
 \int_{B_{t'-t}(0)} & \biggl(
      \displaystyle\int_{B_s}  \frac{|\tau_{h}V_p(Du)|^{2}}{|h|^{2\alpha}} dx
      \biggr)^\frac{\sigma}{2} \frac{dh}{|h|^n}\notag\\
    \le &  
 \theta
  \biggl( \displaystyle\int_{B_{t'}}(1+|D u|)^{p+2\alpha}dx \biggr)^{\frac{\sigma}{2}}
  + C_\theta
  \notag\\
  &
  +\dfrac{C_{\theta}}{(t-s)^{p''}}
  \biggl( \displaystyle\int_{B_t}|D (u-\psi)|^{\frac{p+2\alpha}{p+1+\alpha-q}}dx \biggr)^{\frac{\sigma}{2}}
  +\theta
  \biggl( \displaystyle\int_{B_{t'}}(1+|D u|)^{p+2\alpha}dx \biggr)^{\frac{\sigma}{2}}
\notag\\
&+ \theta  \biggl(  \displaystyle\int_{B_t}(1+|Du|)^{p+2\alpha}dx\biggr)^{\frac{\sigma}{2}}
+ \theta
\left( \ibt(1+|Du|)^{p+2\alpha} dx \right)^\frac{\sigma}{2} \notag\\
&+\frac{C_{\theta}}{(t-s)^{p^*}} 
\left( \ibtt |D (u- \psi)|^\frac{p+2\alpha}{p+1+\alpha-q} dx \right)^\frac{(p+1+\alpha-q)\sigma}{2(p+1+2\alpha-q)} 
+ \theta
\left( \ibt(1+|Du|)^{p+2\alpha} dx \right)^\frac{\sigma}{2}, \label{estimates0.}
\end{align}
for $0 < \theta < 1$, where we denote $p'':= \frac{p+2\alpha}{p+1+\alpha-q}$ and $p^*= \frac{p+2\alpha}{p+1+2\alpha-q}$.

We estimate the second and the penultimate integral appearing in the right hand side of estimate  as follows 

\begin{align}
    & \dfrac{C_\theta}{(t-s)^{p''}} 
    \biggl(
    \displaystyle\int_{B_t}|D (u-\psi)|^{\frac{p+2\alpha}{p+1+\alpha-q}}dx \biggr)^\frac{\sigma}{2}\notag\\
    & \le \dfrac{C_\theta}{(t-s)^{p''}} \biggl(  \displaystyle\int_{B_t}|D u|^{\frac{p+2\alpha}{p+1+\alpha-q}}dx
    \biggr)^\frac{\sigma}{2}
    + \dfrac{C_\theta}{(t-s)^{p''}} 
    \biggl(
    \displaystyle\int_{B_t}|D \psi|^{\frac{p+2\alpha}{p+1+\alpha-q}}dx
    \biggr)^\frac{\sigma}{2}
    \notag\\
    & \le \theta \biggl( \displaystyle\int_{B_t}|D u|^{p+2\alpha}dx
    \biggr)^\frac{\sigma}{2}+ \dfrac{C_\theta(L)}{(t-s)^{\frac{\tilde{p}\sigma}{2}}}   
    |B_R|^\frac{\sigma}{2}
    +  \dfrac{C_\theta}{(t-s)^{p''}} 
    \biggl(
    \displaystyle\int_{B_R}|D \psi|^{\frac{p+2\alpha}{p+1+\alpha-q}}dx
    \biggr)^\frac{\sigma}{2} \label{estimates1.}
\end{align}
and similarly

\begin{align}
    &\frac{C_\theta}{(t-s)^{p^*}} 
\left( \ibtt |D (u- \psi)|^\frac{p+2\alpha}{p+1+\alpha-q} dx \right)^\frac{\sigma(p+1+\alpha-q)}{2(p+1+2\alpha-q)} \notag\\
 & \le  \frac{C_\theta}{(t-s)^{p'''}}
 + C_\theta \biggl(
 \ibtt |D (u- \psi)|^\frac{p+2\alpha}{p+1+\alpha-q} dx
 \biggr)^\frac{\sigma}{2}
 \notag\\
 & \le  
    \frac{C_\theta}{(t-s)^{p'''}}+
    \theta \biggl( \displaystyle\int_{B_{t'}}|D u|^{p+2\alpha}dx
    \biggr)^\frac{\sigma}{2}+
    |B_R|^\frac{\sigma}{2}
    +  C_\theta 
    \biggl(
    \displaystyle\int_{B_{R}}|D \psi|^{\frac{p+2\alpha}{p+1+\alpha-q}}dx
    \biggr)^\frac{\sigma}{2} \label{estimates2.1}
\end{align}
where we set $\tilde{p}= \frac{p+2\alpha}{p+\alpha-q}$ and $p''':= \frac{p+2\alpha}{\alpha}$.\\
Inserting estimates \eqref{estimates1.} and \eqref{estimates2.1} in \eqref{estimates0.}, we obtain

\begin{align}
 \int_{B_{t'-t}(0)} & \biggl(
      \displaystyle\int_{B_s}  \frac{|\tau_{h}V_p(Du)|^{2}}{|h|^{2\alpha}} dx
      \biggr)^\frac{\sigma}{2} \frac{dh}{|h|^n}\notag\\
    \le &  
 4\theta
  \biggl( \displaystyle\int_{B_{t}}(1+|D u|)^{p+2\alpha}dx \biggr)^{\frac{\sigma}{2}}
  + 3 \theta \biggl( \displaystyle\int_{B_{t'}}(1+|D u|)^{p+2\alpha}dx \biggr)^{\frac{\sigma}{2}} \notag\\
 &
  + C_\theta
 + \dfrac{C_\theta(L)}{(t-s)^{\frac{\tilde{p}\sigma}{2}}}   
    |B_R|^\frac{\sigma}{2}
    +  \dfrac{C_\theta}{(t-s)^{p''}} 
    \biggl(
    \displaystyle\int_{B_R}|D \psi|^{\frac{p+2\alpha}{p+1+\alpha-q}}dx
    \biggr)^\frac{\sigma}{2} \notag\\
    &+ \frac{C_\theta}{(t-s)^{p'''}}+
    |B_R|^\frac{\sigma}{2}
    +  C_\theta 
    \biggl(
    \displaystyle\int_{B_{R}}|D \psi|^{\frac{p+2\alpha}{p+1+\alpha-q}}dx
    \biggr)^\frac{\sigma}{2}. \label{tauhvpdu}
\end{align}
Now, by virtue of Proposition \ref{higherintfin}, we infer the following inequality
\begin{align}
    \biggl(\int_{B_{\rho}} |Du|^{p+2\alpha} dx \biggr)^\frac{\sigma}{2}
    \le & 
    C \Vert u \Vert^{\sigma\gamma}_{L^\infty(B_s)}
    \int_{B_{t'-t}(0)} \biggl(\displaystyle\int_{B_s}  \dfrac{|\tau_{h}V_p(Du)|^{2}}{|h|^{2 \alpha}} dx \biggr)^\frac{\sigma}{2} \frac{dh}{|h|^n} 
    \notag\\
    &+ \frac{C}{(s-\rho)^{\sigma p}}
    \Vert u \Vert^{\sigma\gamma}_{L^\infty(B_s)}
    \Vert u \Vert^{\frac{\sigma p}{2}}_{W^{1,p}(B_s)}.\label{interpolazione10}
\end{align}
Combining inequalities \eqref{tauhvpdu} and \eqref{interpolazione10} and arguing as in the proof of Theorem \ref{approximation}, we obtain
\begin{align}
   \biggl( \int_{B_{R/4}} |Du|^{p+2\alpha} dx \biggr)^\frac{\sigma}{2} 
     \leq &
    C
     (\Vert\psi\Vert_{L^\infty(B_R)}+\Vert u \Vert_{W^{1,p}(B_{R})})^\pi \notag\\
     & \cdot  \biggl(
     \Vert \{g_k \}_k \Vert_{l^\sigma(L^\frac{p+2\alpha}{p+\alpha-q}(B_{R}))}^{\sigma} +\Vert D\psi \Vert_{B^\gamma_{\frac{p+2\alpha}{p+1+\alpha-q},\infty}(B_R)} +1
     \biggr)^\pi , 
\end{align}
which yields
\begin{align*}
\int_{B_{t'-t}(0)} \biggl(
\displaystyle\int_{B_{R/4}}\frac{ |\tau_h V_p(Du)|^2}{|h|^{2 \alpha}} dx 
\biggr)^\frac{\sigma}{2} \frac{dh}{|h|^n}
     \le &  C (\Vert\psi\Vert_{L^\infty(B_R)}+\Vert u \Vert_{W^{1,p}(B_{R})})^\pi \notag\\
     & \cdot 
     \biggl(
     \Vert \{g_k \}_k \Vert_{l^\sigma(L^\frac{p+2\alpha}{p+\alpha-q}(B_{R}))}^{\sigma}+\Vert D\psi \Vert_{B^\gamma_{\frac{p+2\alpha}{p+1+\alpha-q},\infty}(B_R)}+1
     \biggr)^\pi,
\end{align*}
for every $t'-t \le \frac{R}{4}$. Hence, we eventually get
\begin{align*}
\int_{B_{R/4}(0)} \biggl(
\displaystyle\int_{B_{R/4}} \frac{|\tau_h V_p(Du)|^2}{|h|^{2 \alpha}} dx 
\biggr)^\frac{\sigma}{2} \frac{dh}{|h|^n}
     \le &  C (\Vert\psi\Vert_{L^\infty(B_R)}+\Vert u \Vert_{W^{1,p}(B_{R})})^\pi \notag\\
     & \cdot 
     \biggl(
     \Vert \{g_k \}_k \Vert_{l^\sigma(L^\frac{p+2\alpha}{p+\alpha-q}(B_{R}))}^{\sigma} +\Vert D\psi \Vert_{B^\gamma_{\frac{p+2\alpha}{p+1+\alpha-q},\infty}(B_R)}+1
     \biggr)^\pi,
\end{align*}
for constants $C := C(n,p,q, \nu, L,R)$ and $\pi:= \pi (n,p,q,\alpha,\sigma)$. 
\qed

\end{document}